\input amstex
\magnification=\magstephalf
\input amsppt.sty

\hsize=15truecm
\vsize=22truecm

\TagsOnRight
\def\lg{\langle}
\def\rg{\rangle}

\def\lra{\longrightarrow}
\def\al{\alpha}
\def\be{\beta}

\def\vn{\varepsilon}

\def\ep{\epsilon}
\def\ot{\otimes}

\def\om{\omega}
\def\lra{\longrightarrow}
\rightheadtext{Drinfeld Doubles and Lusztig Symmetries}
\leftheadtext{N. Bergeron, Y. Gao, N. Hu}

\topmatter

\title\nofrills{DRINFEL'D DOUBLES AND LUSZTIG'S SYMMETRIES OF TWO-PARAMETER QUANTUM GROUPS}
\endtitle
\author Nantel Bergeron$\,^{a,}\,^*$, Yun Gao$\,^{a,}\,^{**}$,  Naihong Hu$\,^{b,}\,^\dagger$ \endauthor

\date First Version on Sept. 1, 2004, Revised Version on April 15, 2005\enddate
\footnote " "{$^*$Supported in part by the NSERC and the CRC of
Canada
 \newline
 \indent $^{**}$Supported in part by the NSERC of Canada
\newline
\indent $^\dagger$Supported by the NNSF (Grant 10431040), the
TRAPOYT and the SFUDP from the MOE of China, the Shanghai Priority
Academic Discipline from the SEC.}

\affil
$^a$Department of Mathematics and Statistics, York University \\
 Toronto, Ontario,  Canada M3J 1P3 \\
$^b$Department of Mathematics, East China Normal University \\
Shanghai 200062, People's Republic of China \\
\endaffil

\abstract  We find the defining structures of two-parameter
quantum groups $U_{r,s}(\frak g)$ corresponding to the orthogonal
and the symplectic Lie algebras, which are realized as Drinfel'd
doubles.  We further investigate the environment conditions upon
which the Lusztig's symmetries exist between $(U_{r,s}(\frak g),
\langle\,,\rangle)$ and its associated object $(U_{s^{-1},
r^{-1}}(\frak g), \langle\,|\,\rangle)$.

\endabstract

\keywords $2$-parameter quantum group, Hopf (skew-)dual pairing,
Drinfel'd double, Lusztig's symmetry
\endkeywords
\subjclass 17B37, 81R50 \endsubjclass

\endtopmatter
\document

\baselineskip13pt

\bigskip
\heading {$\bold 0$. \ Introduction}\endheading
\bigskip\medskip

In the early 1990's, many authors investigated two-parameter or
multiparameter quantum groups. These authors include Kulish
\cite{14}, Reshetikhin \cite{18}, Sudbery \cite{21}, Takeuchi
\cite{22},  Artin, Schelter and  Tate \cite{1}, Du, Parshall and
Wang \cite{9}, Dobrev and Parashar \cite{7}, Jing \cite{11}, Chin
and Musson \cite{6}, etc. (for details, see the introduction of
\cite{3} and references therein). Their works focused on
quantized function algebras and quantum enveloping algebras only
for type $A$ cases. In 2001, Benkart and Witherspoon \cite{3},
motivated by the work on generalizations of algebras generated by
the down and up operators on posets
 (see down-up algebras in  \cite{2}), obtained the
structure of two-parameter quantum enveloping algebras
corresponding to the general linear Lie algebra $\frak {gl}_n$ and
the special linear Lie algebra $\frak{sl}_n$ (which was also
studied earlier by Takeuchi in \cite{22} with a different
motivation).  These two-parameter quantum enveloping algebras of type $A$ were proved to have Drinfel'd doubles
structures. Furthermore, in \cite{4} they studied the
finite-dimensional weight representation theory under the
assumption that $rs^{-1}$ is not a root of unity, and showed that
the finite-dimensional weight modules are completely reducible.
This is analogous to the classical semisimple Lie algebras $\frak
g$ and their one-parameter quantum groups $U_q(\frak g)$ of
Drinfel'd-Jimbo type (with $q$ generic). Since then, a systematic
study for the two-parameter quantum group of type $A$ has been
further developed by Benkart and Witherspoon, etc. (see \cite{5}
and references therein). However, up to now, how to find the
defining relations of two-parameter quantum groups of other types
was still open and an interesting task.

In the present paper, we give the definitions of two-parameter
quantum groups of types $B$, $C$, $D$, which are new to our
knowledge. We also prove that these quantum groups we have constructed
have a Drinfel'd doubles structures. We then investigate the
Lusztig's symmetry properties for these two-parameter quantum
groups for $B$, $C$ and $D$ types as well as for type $A$.
Particularly, we derive some interesting $(r,s)$-identities held
in $U_{r,s}(\frak g)$ (see Lemma 3.6).

The paper is organized as follows. In Section 1, we present  the
definitions of the two-parameter quantum groups corresponding to
the orthogonal Lie algebras $\frak {so}_{2n+1}$ or
$\frak{so}_{2n}$ and the symplectic Lie algebras $\frak
{sp}_{2n}$, together with the Hopf algebra structure. As it  was done
for type $A$ case in \cite{3}, we prove in Section 2 that
two-parameter quantum groups $U_{r,s}(\frak g)$ are characterized
as Drinfel'd doubles $\Cal D(\Cal B, \Cal B')$ of the Hopf
subalgebras $\Cal B$ (upper part) and $\Cal B'$ (lower part) with
respect to a skew-dual paring. As a by-product of the double
structure, one naturally gets a standard triangular decomposition
of $U_{r,s}(\frak g)$ (see Corollary 2.6), which plays a crucial
role in finite dimensional weight representation theory of
$U_{r,s}(\frak g)$. The double structure also leads to a direct
approach to the triangular decomposition of quantum groups
$U_q(\frak g)$ of Drinfel'd-Jimbo type, which, originally, was
obtained nontrivially by Lusztig \cite{16} and Rosso \cite{20} in
different ways. Section 3 is devoted to investigating the
Lusztig's symmetries for the two-parameter quantum groups. A
striking feature of these symmetries is that they exist as $\Bbb
Q$-isomorphisms between $U_{r,s}(\frak g)$ and the associated
object $U_{s^{-1},r^{-1}}(\frak g)$ only when $\text{\rm
rank}\,(\frak g)=2$. When $\text{\rm rank}\,(\frak g)>2$,  the
sufficient and necessary condition for the existence of Lusztig's
symmetries between $U_{r,s}(\frak g)$ and its associated object
forces $U_{r,s}(\frak g)$ to take the ``one-parameter" form
$U_{q,q^{-1}}(\frak g)$ where $r=s^{-1}=q$. In other words, we
prove that when $\text{\rm rank}\,(\frak g)>2$, the Lusztig's
symmetries exist only for the one-parameter quantum groups $U_{q,
q^{-1}}(\frak g)$ as $\Bbb Q(q)$-automorphisms (rather than merely
$\Bbb Q$-isomorphisms). In this case, these symmetries coincide,
modulo identification, with the usual Lusztig symmetries on
quantum groups $U_q(\frak g)$ of Drinfel'd-Jimbo type (see
\cite{12, 13, 15}). Some necessary calculation data are collected
in Section 4, which are useful to the proof of Proposition 2.3.

\bigskip
\bigskip
\heading {$\bold 1$. Two-parameter Quantum Groups of Types
$B,\,C,\,D$}
\endheading
\bigskip

Let $\Bbb K=\Bbb Q(r,s)$ denote a field of rational functions with
two indeterminates $r$, $s$, or a subfield of $\Bbb C$ with
two-parameters $r,\,s$ with assumption $r^2\ne s^2$.

In order to reveal more about the defining structure of the
``two-parameter" quantum groups, we present the definitions case
by case. In particular, the $(r,s)$-Serre relations for the $D_n$
case have some special features for the vertices $n-1$ and $n$,
which degenerate into the usual commutative relations for the
one-parameter case (i.e., when $rs=1$ in $(D5)$ below).

\smallskip
(I) \ Assume $\Psi$ is a finite root system of type $B_n$ with
$\Pi$ a base of simple roots. Regard $\Psi$ as a subset of a
Euclidean space $E=\Bbb R^n$ with an inner product $(\,,)$. Let
$\ep_1, \cdots, \ep_n$ denote an orthonormal basis of $E$, and
suppose $\Pi=\{\al_i=\ep_i-\ep_{i+1}\mid 1\le
i<n\}\cup\{\al_n=\ep_n\}$ and $\Psi=\{\pm\ep_i\pm\ep_j\mid 1\le
i\ne j\le n\}\cup\{\pm\ep_i\mid 1\le i\le n\}$. In this case, set
$r_i=r^{(\al_i, \al_i)}$ and $s_i=s^{(\al_i, \al_i)}$ so that
$r_1=\cdots=r_{n-1}=r^2$, $r_n=r$ and $s_1=\cdots=s_{n-1}=s^2$,
$s_n=s$.

Let $U=U_{r,s}(\frak{so}_{2n+1})$ be the unital associative
algebra over $\Bbb Q(r,s)$ generated by $e_i, f_i,
\omega_i^{\pm 1}, {\om_i'}^{\pm 1}$ $(1\le i\le n)$, subject to
the following relations $(B1)$---$(B7)$:

\smallskip
$(B1)$ \ The $\om_i^{\pm 1}, {\om_j'}^{\pm1}$ all commute with one
another and $\om_i\om_i^{-1}=\om_j'{\om_j'}^{-1}=1$.

\smallskip
$(B2)$ \ For $1\le i\le n$ and $1\le j<n$,  we have
$$\gather
\om_je_i\om_j^{-1}=r_j^{(\ep_j, \al_i)}s_j^{(\ep_{j+1},\al_i)}e_i,
\qquad \om_jf_i\om_j^{-1}=r_j^{-(\ep_j,\al_i)}s_j^{-(\ep_{j+1},\al_i)}f_i,\\
\om_ne_j\om_n^{-1}=r_n^{2(\ep_n, \al_j)}e_j, \qquad
\om_nf_j\om_n^{-1}=r_n^{-2(\ep_n,\al_j)}f_j,\\
\om_ne_n\om_n^{-1}=r_n^{(\ep_n,\al_n)}s_n^{-(\ep_n,\al_n)}e_n,
\qquad
\om_nf_n\om_n^{-1}=r_n^{-(\ep_n,\al_n)}s_n^{(\ep_n,\al_n)}f_n.
\endgather
$$

$(B3)$ \ For $1\le i\le n$ and $1\le j<n$, we have
$$\gather
\om_j'e_i{\om_j'}^{-1}=s_j^{(\ep_j,\al_i)}r_j^{(\ep_{j+1},\al_i)}e_i,
\qquad \om_j'f_i{\om_j'}^{-1}=s_j^{-(\ep_j,\al_i)}r_j^{-(\ep_{j+1},\al_i)}f_i,\\
\om_n'e_j{\om_n'}^{-1}=s_n^{2(\ep_n, \al_j)}e_j, \qquad
\om_n'f_j{\om_n'}^{-1}=s_n^{-2(\ep_n,\al_j)}f_j,\\
\om_n'e_n{\om_n'}^{-1}=s_n^{(\ep_n,
\al_n)}r_n^{-(\ep_n,\al_n)}e_n, \qquad
\om_n'f_n{\om_n'}^{-1}=s_n^{-(\ep_n,
\al_n)}r_n^{(\ep_n,\al_n)}f_n.
\endgather
$$

$(B4)$ \ For $1\le i, j\le n$, we have
$$
[\,e_i, f_j\,]=\delta_{ij}\frac{\om_i-\om_i'}{r_i-s_i}.
$$

\smallskip
$(B5)$ \ For any $i, j$ with $|\,i-j\,|>1$, we have the
$(r,s)$-Serre relations:
$$
[\,e_i, e_j\,]=[\,f_i, f_j\,]=0.
$$

$(B6)$ \ For $1\le i< n$, $1\le j< n-1$, we have the $(r,s)$-Serre
relations:
$$\gather
e_i^2e_{i+1}-(r_i{+}s_i)\,e_ie_{i+1}e_i+(r_is_i)\,e_{i+1}e_i^2=0,\\
e_{j+1}^2e_j-(r_{j+1}^{-1}{+}s_{j+1}^{-1})\,e_{j+1}e_je_{j+1}+(r_{j+1}^{-1}s_{j+1}^{-1})\,e_je_{j+1}^2=0,\\
e_n^3e_{n-1}-(r_n^{-2}{+}r_n^{-1}s_n^{-1}{+}s_n^{-2})\,
e_n^2e_{n-1}e_n+(r_n^{-1}s_n^{-1})(r_n^{-2}{+}r_n^{-1}s_n^{-1}{+}s_n^{-2})\,e_ne_{n-1}e_n^2\\
-(r_n^{-3}s_n^{-3})\,e_{n-1}e_n^3=0.
\endgather
$$

$(B7)$ \ For $1\le i< n$, $1\le j< n-1$, we have the $(r,s)$-Serre
relations:
$$\gather
f_{i+1}f_i^2-(r_i{+}s_i)\,f_if_{i+1}f_i+(r_is_i)\,f_i^2f_{i+1}=0,\\
f_jf_{j+1}^2-(r_{j+1}^{-1}{+}s_{j+1}^{-1})\,f_{j+1}f_jf_{j+1}+(r_{j+1}^{-1}s_{j+1}^{-1})\,f_{j+1}^2f_j=0,\\
f_{n-1}f_n^3-(r_n^{-2}{+}r_n^{-1}s_n^{-1}{+}s_n^{-2})\,f_nf_{n-1}f_n^2+
(r_n^{-1}s_n^{-1})(r_n^{-2}{+}r_n^{-1}s_n^{-1}{+}s_n^{-2})\,f_n^2f_{n-1}f_n\\
-(r_n^{-3}s_n^{-3})\,f_n^3f_{n-1}=0.
\endgather
$$

(II) \ Assume $\Psi$ is a finite root system of type $C_n$ with
$\Pi$ a base of simple roots. Regard $\Psi$ as a subset of a
Euclidean space $E=\Bbb R^n$ with an inner product $(\,,)$. Let
$\ep_1, \cdots, \ep_n$ denote an orthonormal basis of $E$, and
suppose $\Pi=\{\al_i=\ep_i-\ep_{i+1}\mid 1\le
i<n\}\cup\{\al_n=2\ep_n\}$ and $\Psi=\{\pm\ep_i\pm\ep_j\mid 1\le
i\ne j\le n\}\cup\{2\ep_i\mid 1\le i\le n\}$. In this case, set
$r_i=r^{\frac{(\al_i, \al_i)}2}$ and $s_i=s^{\frac{(\al_i, \al_i)}2}$
so that $r_1=\cdots=r_{n-1}=r$, $r_n=r^2$ and
$s_1=\cdots=s_{n-1}=s$, $s_n=s^2$.

Let $U=U_{r,s}(\frak{sp}_{2n})$ be the unital associative algebra
over $\Bbb Q(r,s)$ generated by $e_i$, $f_i$, $\omega_i^{\pm
1}$, ${\om_i'}^{\pm 1}$ $(1\le i\le n)$, subject to the following
relations $(C1)$---$(C7)$:

\smallskip
$(C1)$ \ The $\om_i^{\pm 1}, {\om_j'}^{\pm1}$ all commute with one
another and $\om_i\om_i^{-1}=\om_j'{\om_j'}^{-1}=1$.

\smallskip
$(C2)$ \ For $1\le i\le n$ and $1\le j<n$, we have
$$\gather
\om_je_i\om_j^{-1}=r^{(\ep_j, \al_i)}s^{(\ep_{j+1},\al_i)}e_i,
\qquad \om_jf_i\om_j^{-1}=r^{-(\ep_j,\al_i)}s^{-(\ep_{j+1},\al_i)}f_i,\\
\om_ne_j\om_n^{-1}=r^{2(\ep_n, \al_j)}e_j, \qquad
\om_nf_j\om_n^{-1}=r^{-2(\ep_n,\al_j)}f_j,\\
\om_ne_n\om_n^{-1}=r^{(\ep_n,\al_n)}s^{-(\ep_n,\al_n)}e_n, \qquad
\om_nf_n\om_n^{-1}=r^{-(\ep_n,\al_n)}s^{(\ep_n,\al_n)}f_n.
\endgather
$$

$(C3)$ \ For $1\le i\le n$ and $1\le j<n$, we have
$$\gather
\om_j'e_i{\om_j'}^{-1}=s^{(\ep_j,\al_i)}r^{(\ep_{j+1},\al_i)}e_i,
\qquad \om_j'f_i{\om_j'}^{-1}=s^{-(\ep_j,\al_i)}r^{-(\ep_{j+1},\al_i)}f_i,\\
\om_n'e_j{\om_n'}^{-1}=s^{2(\ep_n, \al_j)}e_j, \qquad
\om_n'f_j{\om_n'}^{-1}=s^{-2(\ep_n,
\al_j)}f_j,\\
\om_n'e_n{\om_n'}^{-1}=s^{(\ep_n, \al_n)}r^{-(\ep_n,\al_n)}e_n,
\qquad \om_n'f_n{\om_n'}^{-1}=s^{-(\ep_n,
\al_n)}r^{(\ep_n,\al_n)}f_n.
\endgather
$$

$(C4)$ \ For $1\le i, j\le n$, we have
$$
[\,e_i, f_j\,]=\delta_{ij}\frac{\om_i-\om_i'}{r_i-s_i}.
$$

\smallskip
$(C5)$ \ For any $i, j$ with $|\,i-j\,|>1$, we have the
$(r,s)$-Serre relations:
$$[\,e_i, e_j\,]=[\,f_i, f_j\,]=0.$$

$(C6)$ \ For $1\le i<n-1$, we have the $(r,s)$-Serre relations:
$$\gather
e_i^2e_{i+1}-(r{+}s)\,e_ie_{i+1}e_i+(rs)\,e_{i+1}e_i^2=0,\\
e_{n-1}^3e_n-(r^2{+}rs{+}s^2)\,e_{n-1}^2e_ne_{n-1}+(rs)\,
(r^2{+}rs{+}s^2)\,e_{n-1}e_ne_{n-1}^2-(rs)^3\,e_ne_{n-1}^3=0,\\
e_{i+1}^2e_i-(r^{-1}{+}s^{-1})\,e_{i+1}e_ie_{i+1}+(r^{-1}s^{-1})\,e_ie_{i+1}^2=0,\\
e_n^2e_{n-1}-(r_n^{-1}{+}s_n^{-1})\,e_ne_{n-1}e_n+(r_n^{-1}s_n^{-1})\,e_{n-1}e_n^2=0.
\endgather
$$

$(C7)$ \ For $1\le i<n-1$, we have the $(r,s)$-Serre relations:
$$\gather
f_{i+1}f_i^2-(r{+}s)\,f_if_{i+1}f_i+(rs)\,f_i^2f_{i+1}=0,\\
f_nf_{n-1}^3-(r^2{+}rs{+}s^2)\,f_{n-1}f_nf_{n-1}^2
+(rs)\,(r^2{+}rs{+}s^2)\,f_{n-1}^2f_nf_{n-1}-(rs)^3\,f_{n-1}^3f_n=0,\\
f_if_{i+1}^2-(r^{-1}{+}s^{-1})\,f_{i+1}f_if_{i+1}+(r^{-1}s^{-1})\,f_{i+1}^2f_i=0,\\
f_{n-1}f_n^2-(r_n^{-1}{+}s_n^{-1})\,f_nf_{n-1}f_n+(r_n^{-1}s_n^{-1})\,f_n^2f_{n-1}=0.
\endgather
$$

(III) \ Assume $\Psi$ is a finite root system of type $D_n$ with
$\Pi$ a base of simple roots. Regard $\Psi$ as a subset of a
Euclidean space $E=\Bbb R^n$ with an inner product $(\,,)$. Let
$\ep_1, \cdots, \ep_n$ denote an orthonormal basis of $E$, and
suppose $\Pi=\{\al_i=\ep_i-\ep_{i+1}\mid 1\le
i<n\}\cup\{\al_n=\ep_{n-1}+\ep_n\}$ and
$\Psi=\{\pm\ep_i\pm\ep_j\mid 1\le i\ne j\le n\}$. In this case,
set $r_i=r^{\frac{(\al_i, \al_i)}2}$ and $s_i=s^{\frac{(\al_i,
\al_i)}2}$ so that $r_1=\cdots=r_{n}=r$ and
$s_1=\cdots=s_{n}=s$.

Let $U=U_{r,s}(\frak{so}_{2n})$ be the unital associative algebra
over $\Bbb Q(r,s)$ generated by $e_i$, $f_i$, $\omega_i^{\pm
1}$, ${\om_i'}^{\pm 1}$ $(1\le i\le n)$, subject to the following
relations $(D1)$---$(D7)$:

\smallskip
$(D1)$ \ The $\om_i^{\pm 1}, {\om_j'}^{\pm1}$ all commute with one
another and $\om_i\om_i^{-1}=\om_j'{\om_j'}^{-1}=1$.

\smallskip
$(D2)$ \ For $1\le i\le n$, $1\le j<n$, and $1\le k\, (\ne
n{-}1)\le n$, we have
$$\gather
\om_je_i\om_j^{-1}=r^{(\ep_j, \al_i)}s^{(\ep_{j+1},\al_i)}e_i,
\qquad \om_jf_i\om_j^{-1}=r^{-(\ep_j,\al_i)}s^{-(\ep_{j+1},\al_i)}f_i,\\
\om_ne_k\om_n^{-1}=r^{(\ep_{n-1}, \al_k)}s^{-(\ep_n, \al_k)}e_k,
\qquad
\om_nf_k\om_n^{-1}=r^{-(\ep_{n-1}, \al_k)}s^{(\ep_n, \al_k)}f_k,\\
\om_ne_{n-1}\om_n^{-1}=r^{(\ep_n,\al_{n-1})}s^{-(\ep_{n-1},\al_{n-1})}e_{n-1},
\qquad
\om_nf_{n-1}\om_n^{-1}=r^{-(\ep_n,\al_{n-1})}s^{(\ep_{n-1},\al_{n-1})}f_{n-1}.
\endgather
$$

$(D3)$ \ For $1\le i\le n$, $1\le j<n$, and $1\le k\, (\ne
n{-}1)\le n$, we have
$$\gather
\om_j'e_i{\om_j'}^{-1}=s^{(\ep_j,\al_i)}r^{(\ep_{j+1},\al_i)}e_i,
\qquad \om_j'f_i{\om_j'}^{-1}=s^{-(\ep_j,\al_i)}r^{-(\ep_{j+1},\al_i)}f_i,\\
\om_n'e_k{\om_n'}^{-1}=s^{(\ep_{n-1}, \al_k)}r^{-(\ep_n,
\al_k)}e_k, \qquad \om_n'f_k{\om_n'}^{-1}=s^{-(\ep_{n-1},
\al_k)}r^{(\ep_n, \al_k)}f_k,\\
\om_n'e_{n-1}{\om_n'}^{{-}1}=s^{(\ep_n,\al_{n{-}1})}r^{{-}(\ep_{n{-}1},\al_{n{-}1})}e_{n{-}1},
\,\quad
\om_n'f_{n-1}{\om_n'}^{{-}1}=s^{{-}(\ep_n,\al_{n{-}1})}r^{(\ep_{n{-}1},\al_{n{-}1})}f_{n{-}1}.
\endgather
$$

$(D4)$ \ For $1\le i, j\le n$, we have
$$
[\,e_i, f_j\,]=\delta_{ij}\frac{\om_i-\om_i'}{r-s}.
$$

\smallskip
$(D5)$ \ For any $1\le i\ne j\le n$ but $(i, j)\not\in\{(n{-}1,
n), (n, n{-}1)\}$ with $a_{ij}=0$, we have the $(r,s)$-Serre
relations:
$$\gather
[\,e_i, e_j\,]=[\,f_i, f_j\,]=0,\\
e_{n-1}e_n=rs\,e_ne_{n-1},\qquad f_nf_{n-1}=rs\,f_{n-1}f_n.
\endgather
$$

$(D6)$ \ For $1\le i<j\le n$ with $a_{ij}=-1$, we have the
$(r,s)$-Serre relations:
$$\gather
e_i^2e_j-(r{+}s)\,e_ie_je_i+(rs)\,e_je_i^2=0,\\
e_j^2e_i-(r^{-1}{+}s^{-1})\,e_je_ie_j+(r^{-1}s^{-1})\,e_ie_j^2=0.
\endgather
$$

$(D7)$ \ For $1\le i<j\le n$ with $a_{ij}=-1$, we have the
$(r,s)$-Serre relations:
$$\gather
f_jf_i^2-(r{+}s)\,f_if_jf_i+(rs)\,f_i^2f_j=0,\\
f_if_j^2-(r^{-1}{+}s^{-1})\,f_jf_if_j+(r^{-1}s^{-1})\,f_j^2f_i=0.
\endgather
$$

\medskip
In summary, let $U=U_{r,s}(\frak g)$ for $\frak
g=\frak{so}_{2n+1}$, $\frak{so}_{2n}$ and $\frak{sp}_{2n}$ denote
the two-parameter quantum orthogonal groups and quantum symplectic
groups, respectively.

The following fact is straightforward to check.
\proclaim{Proposition 1.1} \ The algebra $U_{r, s}(\frak g)$
$($\,$\frak g=\frak{so}_{2n+1},\,\frak{so}_{2n}$, or
$\frak{sp}_{2n}$\,$)$ is a Hopf algebra under the
comultiplication, the counit and the antipode below:
$$
\gather \Delta(\om_i^{\pm1})=\om_i^{\pm1}\ot\om_i^{\pm1}, \qquad
\Delta({\om_i'}^{\pm1})={\om_i'}^{\pm1}\ot{\om_i'}^{\pm1},\\
\Delta(e_i)=e_i\ot 1+\om_i\ot e_i, \qquad \Delta(f_i)=1\ot
f_i+f_i\ot \om_i',\\
\vn(\om_i^{\pm1})=\vn({\om_i'}^{\pm1})=1, \qquad
\vn(e_i)=\vn(f_i)=0,\\
S(\om_i^{\pm1})=\om_i^{\mp1}, \qquad
S({\om_i'}^{\pm1})={\om_i'}^{\mp1},\\
S(e_i)=-\om_i^{-1}e_i,\qquad S(f_i)=-f_i\,{\om_i'}^{-1}.
\endgather
$$
\endproclaim

\noindent {\bf Remarks.} \ (1) \ When $r=q$ and $s=q^{-1}$, the
Hopf algebra $U_{q, q^{-1}}(\frak g)$ modulo the Hopf ideal
generated by the elements $\om_i'-\om_i^{-1}$ $(1\le i\le n)$, is
 the one-parameter quantum groups $U_q(\frak g)$ of
Drinfel'd-Jimbo type.

\smallskip
(2) \ As usual, we define respectively the left-adjoint and the right-adjoint
action in the Hopf algebra $U_{r, s}(\frak g)$ as
$$
\text{ad}_{ l}\,a\,(b)=\sum_{(a)}a_{(1)}\,b\,S(a_{(2)}), \qquad
\text{ad}_{ r}\,a\,(b)=\sum_{(a)}S(a_{(1)})\,b\,a_{(2)},
$$
where $\Delta(a)=\sum_{(a)}a_{(1)}\ot a_{(2)}$, for any $a$, $b\in
U_{r,s}(\frak g)$.

Using adjoint actions, the $(r,s)$-Serre relations $(X5)$, $(X6)$
and $(X7)$ (here $X=B$, $C$, $D$) in $U_{r,s}(\frak g)$ become
simply:
$$\gather
\bigl(\text{ad}_l\,e_i\bigr)^{1-a_{ij}}\,(e_j)=0,
\qquad\text{\it for any } \ i\ne j,\tag{X6}\\
\bigl(\text{ad}_r\,f_i\bigr)^{1-a_{ij}}\,(f_j)=0, \qquad\text{\it
for any } \ i\ne j.\tag{X7}
\endgather
$$

According to the data on the prime root systems of the classical
simple Lie algebras, the following basic lemma is clear. This will
play a crucial role in ensuring the compatibility of the defining
relations $(X2)$ \& $(X3)$ for $U_{r,s}(\frak g)$, especially, in
the proof of Theorem 2.5.
\proclaim {Lemma 1.2} \ For the prime
root systems of the Lie algebras $\frak g=\frak {sl}_n$,
$\frak{so}_{2n+1}$, $\frak{so}_{2n}$, and $\frak{sp}_{2n}$, the following identities hold:
$$
\gather (\ep_{j+1},\al_i)=-(\ep_i,\al_j),\, \qquad (i, j<n),
\;\qquad\qquad\quad
\text{\it for any } \ \frak g,\\
(\ep_{j+1},\al_n)=\cases-(\ep_n,\al_j), \quad &\quad (j<n),
\qquad\quad\quad \qquad \text{\it
for } \ \frak g=\frak{so}_{2n+1}, \\
-2(\ep_n,\al_j), \quad &\quad (j<n), \quad \qquad\qquad\quad
\text{\it for } \ \frak g=\frak{sp}_{2n},
\endcases\\
(\ep_j,\al_n)=\cases-(\ep_n,\al_{j-1}),  &\quad (j\le n,\, j\ne n-1),  \\
(\ep_{n-1},\al_{n-1}), &\quad (j=n-1)
\endcases
\,\quad\text{\it for } \ \frak g=\frak{so}_{2n}.
\endgather
$$
\endproclaim

\bigskip\medskip
\heading{$\bold 2$. \ Drinfel'd Quantum Doubles and Rosso
Form}\endheading
\bigskip

\proclaim{Definition 2.1} \ A {\it $($Hopf $)$ dual pairing} of
two Hopf algebras $\Cal A$ and $\Cal U$ $($see \cite{12, 3.2.1},
or \cite{13, 1.2.5}$)$ is a bilinear form $\langle\,,\rangle:\;
\Cal U\times \Cal A\lra \Bbb K$ such that
$$
\gather \langle f, 1_{\Cal A}\rangle=\vn_{\Cal U}(f),\qquad
\langle 1_{\Cal U}, a\rangle=\vn_{\Cal A}(a),\tag{1}\\
\langle f, a_1a_2\rangle=\langle \Delta_{\Cal U}(f), a_1\ot
a_2\rangle, \qquad \langle f_1f_2, a \rangle=\langle f_1\ot f_2,
\Delta_{\Cal A}(a)\rangle,\tag{2}
\endgather
$$
for all $f,\, f_1,\, f_2\in\Cal U$, and $a,\,a_1,\,a_2\in\Cal A$.
Here $\vn_{\Cal U}$ and $\vn_{\Cal A}$ denote the counits of
$\Cal U$ and $\Cal A$, respectively, and $\Delta_{\Cal U}$ and
$\Delta_{\Cal A}$ are their comultiplications.
\endproclaim
A direct consequence of the defining properties above is that the dual
pairing satisfies
$$
\langle S_{\Cal U}(f),a\rangle=\langle f, S_{\Cal A}(a)\rangle,
\quad f\in\Cal U, \ a\in\Cal A,
$$
where $S_{\Cal U},\, S_{\Cal A}$ denote the respective antipodes
of \,$\Cal U$ and $\Cal A$.

\proclaim{Definition 2.2} \ A bilinear form $\langle\,,\rangle:\,
\Cal U\times \Cal A\lra \Bbb K$ is called a skew-dual pairing of
two Hopf algebras $\Cal A$ and $\Cal U$ $($see \cite{13, 8.2.1}$)$
if $\langle\,,\rangle:\, \Cal U^{\text{cop}}\times \Cal A\lra \Bbb
K$ is a dual pairing of the two Hopf algebras $\Cal A$ and $\Cal
U^{\text{cop}}$. Here $\Cal U^{\text{cop}}$ is the Hopf algebra
having the opposite comultiplication to the Hopf algebra $\Cal U$
and $S_{\Cal U^{\text{cop}}}=S^{-1}_{\Cal U}$ if $S_{\Cal U}$ is
invertible.
\endproclaim

Let $\Cal B=B(\frak g)$ (resp. $\Cal B'=B'(\frak g)$\,) denote the
Hopf subalgebra of $U=U_{r,s}(\frak g)$  generated by $e_j$,
$\om_j^{\pm1}$ (resp. $f_j$, ${\om_j'}^{\pm1}$\,) with $1\le j< n$
for $\frak g=\frak{sl}_n$, and with $1\le j\le n$ for $\frak
g=\frak {so}_{2n+1}$, $\frak {so}_{2n}$, or $\frak{sp}_{2n}$. Note
that the case $\frak g=\frak {sl}_n$ was done in \cite{3}. The
main ideas of this section comes from this work. Since our
statement is based on Definition 2.2 and slightly different, we
also include the case $\frak g=\frak{sl}_n$ in Proposition 2.3 and
Theorem 2.5 below.

\proclaim{Proposition 2.3} \ There exists a unique skew-dual
pairing $\langle\,,\rangle:\, B'(\frak g)\times B(\frak g)\lra\Bbb
Q(r,s)$ of the Hopf subalgebras $B(\frak g)$ and $B'(\frak g)$,
for $\frak g=\frak {sl}_n$, $\frak {so}_{2n+1}$, $\frak
{so}_{2n}$, or $\frak{sp}_{2n}$ such that
$$
\gather
\lg f_i, e_j\rg=\delta_{ij}\frac1{s_i-r_i},\qquad \text{\it for any } \ \frak g,\tag{3}\\
\lg \om_i', \om_j\rg= r^{(\ep_j,\al_i)}s^{(\ep_{j+1},\al_i)},
\qquad\qquad\  i< n, \ j<n,  \
\qquad\text{\it for}\quad \frak{sl}_n,\tag{4A}\\
\lg \om_i', \om_j\rg=\cases
r^{2(\ep_j,\al_i)}s^{2(\ep_{j+1},\al_i)}, & \quad\qquad\ i\le n, \
j<n,\cr r^{2(\ep_n, \al_i)}, & \quad\qquad\ i<n, \ j=n,\qquad\
\text{\it for}\quad \frak{so}_{2n+1},\cr r^{(\ep_n,
\al_n)}s^{-(\ep_n, \al_n)}, & \qquad\quad\  i=j=n.
\endcases\tag{4B}\\
\lg \om_i', \om_j\rg=\cases
r^{(\ep_j,\al_i)}s^{(\ep_{j+1},\al_i)}, & \quad\qquad  i\le n, \
j<n,\cr r^{2(\ep_n, \al_i)}, & \quad\qquad i<n, \ j=n, \
\qquad\;\text{\it for}\quad \frak{sp}_{2n},\cr r^{(\ep_n,
\al_n)}s^{-(\ep_n, \al_n)}, & \quad\qquad i=j=n.
\endcases\tag{4C}\\
\lg \om_i', \om_j\rg=\cases
r^{(\ep_j,\al_i)}s^{(\ep_{j+1},\al_i)}, & \ \;  i\le n, \ j<n,\cr
r^{(\ep_{n-1}, \al_i)}s^{-(\ep_n,\al_i)}, & \ \; i\ne n-1, \
j=n,\quad\text{\it for}\quad \frak{so}_{2n},\cr r^{(\ep_n,
\al_{n-1})}s^{-(\ep_{n-1}, \al_{n-1})}, & \ \;i=n-1, j=n.
\endcases\tag{4D}\\
\lg {\om_i'}^{\pm1},\om_j^{-1}\rg=\lg
{\om_i'}^{\pm1},\om_j\rg^{-1}=\lg \om_i',\om_j\rg^{\mp1},\qquad
\text{\it for any } \ \frak g,  \tag{5}
\endgather
$$
and all other pairs of generators are $0$. Moreover, we have $\lg
S(a), S(b)\rg=\lg a, b\rg$ for $a\in\Cal B',\,b\in\Cal B$.
\endproclaim
\demo{Proof} \ The uniqueness assertion is clear, since any
skew-dual pairing of bialgebras is determined by the values on the
generators. We proceed to prove the existence of the pairing.

We begin by defining on the generators a bilinear form $\langle\,,\rangle: {\Cal
B'}^{\text{cop}}\times \Cal B\lra \Bbb Q(r,s)$  satisfying (3), (4X) and (5).
We then extend it to a
bilinear form on ${\Cal B'}^{\text{cop}}\times \Cal B$ by
requiring that (1) and (2) hold for $\Delta_{{\Cal
B'}^{\text{cop}}}=\Delta_{\Cal B'}^{\text{op}}$. We  verify
that the relations in $\Cal B$ and $\Cal B'$ are preserved,
ensuring that the form is well-defined and is a dual pairing of
$\Cal B$ and ${\Cal B'}^{\text{cop}}$.

It is straightforward to check that the bilinear form preserves
all the relations among the $\om_i^{\pm1}$ in $\Cal B$ and the
${\om_i'}^{\pm1}$ in $\Cal B'$. We notice that Lemma 1.2 ensures
the compatibility of the form defined above with the relations of
(X2) and (X3) in $\Cal B$ and $\Cal B'$ respectively. This follows
by definition (from (3), (4X), \& (5)). We are left  to verify
that the form preserves the $(r,s)$-Serre relations in $\Cal B$
and $\Cal B'$. It suffices to show that the form on ${\Cal
B'}^{\text{cop}}\times \Cal B$ preserves the $(r,s)$-Serre
relations in $\Cal B$, because the verification for ${\Cal
B'}^{\text{cop}}$ is similar. To this end, we observe that the
relations involving indices $1\le i,\,j<n$ belong to $\Cal B\cap
U_{r,s}(\frak {sl}_n)=B(\frak {sl}_n)$. This has been checked by
Benkart and Witherspoon in the type $A$ case (see \cite{3, Lemma
2.2}). We have thus reduced the proof to the rank $2$ type $B_2$
and $C_2$, and the rank $4$ type $D_4$.

\medskip
(I) \ Cases $\frak g=\frak {sp}_{2n}$ and $\frak{so}_{2n+1}$: we only need
 to consider the case $n=2$ ($B_2$ and $C_2$).

First, let us show that the form preserves the
$(r,s)$-Serre relation of degree $3$ in $\Cal B$:
$$\gather
\lg\, X,\,
e_1^3e_2-(r^2{+}rs{+}s^2)\,e_1^2e_2e_1+(rs)(r^2{+}rs{+}s^2)\,
e_1e_2e_1^2-(rs)^3e_2e_1^3\,\rg=0,\tag{\text{$C_2$}}\\
\lg\, X,\,
e_2^3e_1-(r^{{-}2}{+}r^{{-}1}s^{{-}1}{+}s^{{-}2})\,e_2^2e_1e_2+
(r^{{-}1}s^{{-}1})(r^{{-}2}{+}r^{{-}1}s^{{-}1}{+}s^{{-}2})\,e_2e_1e_2^2\\
-(r^{{-}1}s^{{-}1})^3e_1e_2^3\,\rg=0,\tag{\text{$B_2$}}
\endgather
$$
where $X$ is any word in the generators of $\Cal B'$. By
definition, the left-hand side of each of the above identities
respectively equals
$$\gather
\lg\, \Delta^{(3)}(X),\, e_1\ot e_1\ot e_1\ot
e_2-(r^2{+}rs{+}s^2)\,e_1\ot e_1\ot e_2\ot e_1\\
+(rs)(r^2{+}rs{+}s^2)\,e_1\ot e_2\ot e_1\ot e_1-(rs)^3e_2\ot
e_1\ot e_1\ot e_1\,\rg,\tag{6C} \\
\lg\, \Delta^{(3)}(X),\, e_2\ot e_2\ot e_2\ot
e_1-(r^{{-}2}{+}r^{{-}1}s^{{-}1}{+}s^{{-}2})\,e_2\ot e_2\ot e_1\ot e_2\\
+(r^{{-}1}s^{{-}1})(r^{{-}2}{+}r^{{-}1}s^{{-}1}{+}s^{{-}2})\,e_2\ot
e_1\ot e_2\ot e_2-(r^{{-}3}s^{{-}3})\,e_1\ot e_2\ot e_2\ot
e_2\,\rg,\tag{6B}
\endgather
$$
where the $\Delta$ corresponds to $\Delta_{\Cal B'}^{\text{op}}$. In order for
any one of these terms to be nonzero, $X$ must involve exactly
three $f_1$ factors, one $f_2$ factor, and arbitrarily many
${\om_j'}^{\pm1}$ factors $(j=1, 2$).

We first consider the following four key cases:

\smallskip
(i) \ For type $C_2$ and $X=f_1^3f_2$, we have
$$\split
\Delta^{(3)}(X)=\,&\bigl(\om_1'\ot\om_1'\ot\om_1'\ot
f_1+\om_1'\ot\om_1'\ot f_1\ot
1+\om_1'\ot f_1\ot 1\ot 1+f_1\ot 1\ot 1\ot 1\bigr)^3\\
&\cdot\bigl(\om_2'\ot\om_2'\ot\om_2'\ot f_2+\om_2'\ot\om_2'\ot
f_2\ot 1+\om_2'\ot f_2\ot 1\ot 1+f_2\ot 1\ot 1\ot 1\bigr).
\endsplit
$$
In Appendix (4.1) we have listed only the terms in the expansion
of $\Delta^{(3)}(X)$ that
 have a non-zero contribution in (6C). Consequently, by properties
(1) \& (2) of $\lg\,,\rg$ on ${\Cal B'}^{\text{cop}}\times \Cal
B$, the pairing (6C) becomes
$$\split
&\lg f_1,e_1\rg^3\lg
f_2,e_2\rg\,\bigl(1+2\lg\om_1',\om_1\rg+2\lg\om_1',\om_1\rg^2+
\lg\om_1',\om_1\rg^3\bigr)\,\cdot\\
&\quad\cdot\bigl(1-r^2\lg\om_1',\om_2\rg\bigr)\bigl(1-rs\lg\om_1',\om_2\rg\bigr)
\bigl(1-s^2\lg\om_1',\om_2\rg\bigr)\\
&=0, \quad(\text{\it since }\lg\om_1',\om_2\rg=r^{-2}).
\endsplit
$$

For type $B_2$, when $X=f_2^3f_1$, we only need to interchange the indices $1$
and $2$, and substitute $r,\,s$ by $r^{-1},\,s^{-1}$ in the
above identity, to obtain the pairing (6B) as
$$\split
&\lg f_1,e_1\rg\lg f_2,e_2\rg^3\,\bigl(1+2\lg\om_2',\om_2\rg+
2\lg\om_2',\om_2\rg^2+\lg\om_2',\om_2\rg^3\bigr)\\
&\quad\cdot\bigl(1-r^{-2}\lg\om_2',\om_1\rg\bigr)\bigl(1-r^{-1}s^{-1}
\lg\om_2',\om_1\rg\bigr)\bigl(1-s^{-2}\lg\om_2',\om_1\rg\bigr)\\
&=0, \quad(\text{\it since }\lg\om_2',\om_1\rg=s^2).
\endsplit
$$

(ii) \ For type $C_2$ and $X=f_2f_1^3$, the relevant terms of
$\Delta^{(3)}(X)$ for (6C) are listed in Appendix (4.2). The
result of (6C) is then
$$\split
&\lg f_1,e_1\rg^3\lg
f_2,e_2\rg\,\bigl(1+2\lg\om_1',\om_1\rg+2\lg\om_1',
\om_1\rg^2+\lg\om_1',\om_1\rg^3\bigr)\,\cdot\\
&\quad\cdot\bigl(\lg\om_2',\om_1\rg-r^2\bigr)\bigl(\lg\om_2',\om_1\rg-rs\bigr)
\bigl(\lg\om_2',\om_1\rg-s^2\bigr)\\
&=0, \quad(\text{\it since }\lg\om_2',\om_1\rg=s^2).
\endsplit
$$

Again for type $B_2$, when $X=f_1f_2^3$ we need to interchange the indices $1$
and $2$, and substitute $r,\,s$ by $r^{-1},\,s^{-1}$ in the
above identity, to obtain the pairing (6B) as
$$\split
&\lg f_1,e_1\rg\lg
f_2,e_2\rg^3\,\bigl(1+2\lg\om_2',\om_2\rg+
2\lg\om_2',\om_2\rg^2+\lg\om_2',\om_2\rg^3\bigr)\\
&\quad\cdot\bigl(\lg\om_1',\om_2\rg-r^{-2}\bigr)
\bigl(\lg\om_1',\om_2\rg-r^{-1}s^{-1}\bigr)\bigl(\lg\om_1',\om_2\rg
-s^{-2}\bigr)\\
&=0, \quad(\text{\it since }\lg\om_1',\om_2\rg=r^{-2}).
\endsplit
$$

(iii) \ For type $C_2$ and $X=f_1^2f_2f_1$,  the relevant terms of
$\Delta^{(3)}(X)$ for (6C) are listed in Appendix (4.3), and (6C)
becomes
$$\split
&\lg f_1,e_1\rg^3\lg
f_2,e_2\rg\,\left\{\bigl(\,1+2\lg\om_1',\om_1\rg
+2\lg\om_1',\om_1\rg^2+\lg\om_1',\om_1\rg^3\,\bigr)
\bigl(\,\lg\om_2',\om_1\rg-(rs)^3\lg\om_1',\om_2\rg^2\,\bigr)\right.\\
&\quad+(r^2+rs+s^2)\,\bigl[\,(rs)\,\lg\om_1',\om_2\rg\,
\bigl(\,1+2\lg\om_1',\om_1\rg+\lg\om_1',\om_1\rg^2\\
&\qquad\qquad\qquad\qquad\qquad
+\lg\om_1',\om_1\rg^2\lg\om_1',\om_2\rg\lg\om_2',\om_1\rg
+\lg\om_1',\om_1\rg^3\lg\om_1',\om_2\rg\lg\om_2',\om_1\rg\,\bigr)\\
&\qquad\qquad\qquad\qquad
-\bigl(\,1+\lg\om_1',\om_1\rg+\lg\om_1',\om_1\rg\lg\om_1',\om_2\rg\lg\om_2',\om_1\rg\\
&\qquad\qquad\qquad\qquad\qquad\left.+2\lg\om_1',\om_1\rg^2\lg\om_1',\om_2\rg\lg\om_2',\om_1\rg
+\lg\om_1',\om_1\rg^3\lg\om_1',\om_2\rg\lg\om_2',\om_1\rg\,\bigr)\,\bigr]\right\}\\
&=0, \quad (\text{\it since } \lg\om_1',\om_1\rg=rs^{-1}, \
\lg\om_1',\om_2\rg=r^{-2}, \ \lg\om_2',\om_1\rg=s^2).
\endsplit
$$

Similarly, for $B_2$ and $X=f_2^2f_1f_2$, interchanging indices
$1$ and $2$, and substituting $r,\,s$ by $r^{-1},\,s^{-1}$ in the
above identity, we obtain that the pairing (6B) is zero.

(iv) \ For type $C_2$ and $X=f_1f_2f_1^2$, the relevant terms of
$\Delta^{(3)}(X)$ for (6C) are listed in Appendix (4.4), and (6C)
becomes
$$
\split &\lg f_1,e_1\rg^3\lg
f_2,e_2\rg\,\bigl(\,1+2\lg\om_1',\om_1\rg+2\lg\om_1',\om_1\rg^2+\lg\om_1',\om_1\rg^3\,\bigr)
\bigl(\,\lg\om_2',\om_1\rg^2-(rs)^3\lg\om_1',\om_2\rg\,\bigr)\\
&\quad+(r^2+rs+s^2)\,\bigl[\,(rs)\,\bigl(\,1+\lg\om_1',\om_1\rg+\lg\om_1',
\om_1\rg\lg\om_1',\om_2\rg\lg\om_2',\om_1\rg\\
&\qquad\qquad\qquad\qquad\qquad
+2\lg\om_1',\om_1\rg^2\lg\om_1',\om_2\rg\lg\om_2',\om_1\rg
+\lg\om_1',\om_1\rg^3\lg\om_1',\om_2\rg\lg\om_2',\om_1\rg\,\bigr)\\
&\qquad\qquad\qquad\qquad
-\lg\om_2',\om_1\rg\,\bigl(\,1+2\lg\om_1',\om_1\rg+\lg\om_1',\om_1\rg^2\\
&\qquad\qquad\qquad\qquad\qquad+\lg\om_1',\om_1\rg^2\lg\om_1',\om_2\rg\lg\om_2',\om_1\rg
+\lg\om_1',\om_1\rg^3\lg\om_1',\om_2\rg\lg\om_2',\om_1\rg\,\bigr)\,\bigr]\\
&=0, \quad (\text{\it since } \lg\om_1',\om_1\rg=rs^{-1}, \
\lg\om_1',\om_2\rg=r^{-2}, \ \lg\om_2',\om_1\rg=s^2).
\endsplit
$$

Again for type $B_2$ and $X=f_2f_1f_2^2$, interchanging indices
$1$ and $2$, and substituting $r,\,s$ by $r^{-1},\,s^{-1}$ in the
above identity, we obtain that the pairing (6B) is also zero.

Finally, we note the fact that if
$X={\om_2'}^pf_1^3{\om_1'}^qf_2{\om_j'}^s$ (resp.
${\om_1'}^qf_2^3{\om_2'}^pf_1{\om_j'}^s$), then (6C) (resp. (6B))
is $\lg\om_2',\om_1\rg^p\lg\om_1',\om_2\rg^q$ times the
corresponding quantity for $X=f_1^3f_2$ for $C_2$ (resp.
$X=f_2^3f_1$ for $B_2$). That means, if $X$ is any word involving
exactly three $f_1$ factors, one $f_2$ factor, and arbitrarily
many factors of ${\om_j'}^{\pm1}$ ($j=1, 2$), then (6C) (resp.
(6B)) is just a scalar multiple of one of the quantities we
have already calculated, and therefore  equal $0$.

\medskip
Next, we verify that the $(r,s)$-Serre relation of degree $2$
in $\Cal B$ are preserved by the form. We do so by showing the identities
$$\gather
\lg X,\,
e_2^2e_1-(r^{-2}{+}s^{-2})\,e_2e_1e_2+(r^{-2}s^{-2})\,e_1e_2^2\rg
=0,\tag{\text{7C}}\\
\lg X, \, e_1^2e_2-(r^2+s^2)\,e_1e_2e_1+(r^2s^2)\,e_2e_1^2\rg
=0,\tag{\text{7B}}\\
\endgather
$$
where $X$ is any word in the generators of $\Cal B'$. It is enough
to consider three monomials: $X=f_2^2f_1$, $f_2f_1f_2$ and
$f_1f_2^2$ for $C_2$ (resp. $X=f_1^2f_2$, $f_1f_2f_1$ and
$f_2f_1^2$ for $B_2$).

\medskip (i) \ For type $C_2$ and $X=f_2^2f_1$, the relevant terms of
$\Delta^{(2)}(X)$ are listed in Appendix (4.5), and the left-hand
side of (7C) equals
$$
\lg f_2, e_2\rg^2\lg
f_1,e_1\rg\,\bigl(1+\lg\om_2',\om_2\rg\,\bigr)
\bigl(1-r^{-2}\lg\om_2',\om_1\rg\,\bigr)\bigl(1-s^{-2}\lg\om_2',\om_1\rg\,\bigr),
$$
which is $0$, since $\lg\om_2',\om_1\rg=s^2$.

For $B_2$ and $X=f_1^2f_2$, the left-hand side of (7B) is
$$
\lg f_1,e_1\rg^2\lg f_2,e_2\rg\,\bigl(1+\lg\om_1',\om_1\rg\,\bigr)
\bigl(1-r^2\lg\om_1',\om_2\rg\,\bigr)\bigl(1-s^2\lg\om_1',\om_2\rg\,\bigr),
$$
which is $0$, since $\lg\om_1',\om_2\rg=r^{-2}$.

\medskip
(ii) \ Similarly, for type $C_2$ and $X=f_1f_2^2$ (resp. type
$B_2$ and $X= f_2f_1^2$), the left-hand side of (7C) and (7B)
respectively are equal to
$$\gather
\lg f_1,e_1\rg\,\lg f_2,e_2\rg^2\bigl(1+\lg\om_2',\om_2\rg\,\bigr)
\bigl(\lg\om_1',\om_2\rg-r^{-2}\bigr)\bigl(\lg\om_1',\om_2\rg-s^{-2}\bigr)=0,
\quad(\text{\it as } \lg\om_1',\om_2\rg=r^{-2}); \\
\lg f_1,e_1\rg^2\lg f_2,e_2\rg\,\bigl(1+\lg\om_1',\om_1\rg\,\bigr)
\bigl(\lg\om_2',\om_1\rg-r^2\bigr)\bigl(\lg\om_2',\om_1\rg-s^2\bigr)=0,
\quad(\text{\it as } \lg\om_2',\om_1\rg=s^2).
\endgather
$$

(iii) \ For type $C_2$ and $X=f_2f_1f_2$, the relevant terms of
$\Delta^{(2)}(X)$ are listed in Appendix (4.6), and the left-hand
side of (7C) equals
$$\split
\lg f_2&,e_2\rg^2\lg
f_1,e_1\rg\bigl\{\,\bigl(\,1+\lg\om_2',\om_2\rg\,\bigr)
\bigl(\,\lg\om_1',\om_2\rg+r^{-2}s^{-2}\lg\om_2',\om_1\rg\,\bigr)\\
&-(r^{-2}+s^{-2})\bigl(\,1+\lg\om_2',\om_2\rg\lg\om_1',\om_2\rg\lg\om_2',\om_1\rg\,\bigr)\,\bigr\},
\endsplit
$$
which is equal to $0$, since $\lg\om_2',\om_2\rg=r^2s^{-2}$,
$\lg\om_1',\om_2\rg=r^{-2}$ and $\lg\om_2',\om_1\rg=s^2$.

Similarly, for $B_2$ and $X=f_1f_2f_1$, the left-hand
side of (7B) is
$$\split
\lg f_1&,e_1\rg^2\lg
f_2,e_2\rg\bigl\{\,\bigl(\,1+\lg\om_1',\om_1\rg\,\bigr)
\bigl(\,\lg\om_2',\om_1\rg+r^2s^2\lg\om_1',\om_2\rg\,\bigr)\\
&-(r^2+s^2)\bigl(\,1+\lg\om_1',\om_1\rg\lg\om_2',
\om_1\rg\lg\om_1',\om_2\rg\,\bigr)\,\bigr\},
\endsplit
$$
which is equal to $0$, since $\lg\om_1',\om_1\rg=r^2s^{-2}$,
$\lg\om_1',\om_2\rg=r^{-2}$ and $\lg\om_2',\om_1\rg=s^2$.

\medskip
A similar process shows that the relations in ${\Cal
B'}^{\text{cop}}$ are preserved in the case when $\frak
g=\frak{sp}_{2n}$, or $\frak{so}_{2n+1}$.

\medskip
(II) \ Case $\frak g=\frak{so}_{2n}$: we only need to consider the
type $D_4$.

For $j=3$ or $4$, we let $\Cal B_j$ (resp. $\Cal B'_j$) denote the subalgebra
generated by $\om_i,\,e_i$ (resp.  $\om_i',\,f_i$) with
$i\in \{1,\,2,\,j\}$. Note that $\Cal B_3\cong\Cal
B_4\cong B(\frak{sl}_4)$ and $\Cal B_3\cap\Cal B_4=B(\frak{sl}_3)\subset
B(\frak{sl}_4)$;  $\Cal B_3'\cong\Cal B_4'\cong
B'(\frak{sl}_4)$ and $\Cal B_3'\cap\Cal
B_4'=B'(\frak{sl}_3)\subset B'(\frak{sl}_4)$.
 The existence of
pairings on ${\Cal B_j'}^{\text{cop}}\times\Cal B_j$ ($j=3,\,4$),
whose restrictions on ${B'(\frak{sl}_3)}^{\text{cop}}\times
B(\frak{sl}_3)$ are consistent, has been proved in \cite{3, Lemma
2.2}. Therefore, in order to show that the form defined on ${\Cal
B'}^{\text{cop}}\times\Cal B$ preserves the $(r,s)$-Serre
relations in $\Cal B$, it remains to verify that the $(r,s)$-Serre
relation in $\Cal B$ involving indices $3,\,4$ in (D5) is
preserved. We need only to consider two cases:
$X=f_if_j$ for $(i,j)=(3,4)$, or $(4,3)$.
$$\split
\lg \,X,& \,e_3e_4-(rs)\,e_4e_3\,\rg\\
&=\lg f_i\om_j'\ot f_j+\om_i'f_j\ot f_i, \, e_3\ot e_4-(rs)\,e_4\ot e_3\rg\\
&=\lg f_i, e_3\rg\lg f_j,e_4\rg-(rs)\,\lg f_i, e_4\rg\lg f_j,
e_3\rg\\
&\quad+\lg\om_i',\om_3\rg\lg f_j,e_3\rg\lg f_i,
e_4\rg-(rs)\,\lg\om_i',\om_4\rg\lg f_j,e_4\rg\lg f_i,e_3\rg\\
&=\cases\frac{1}{(s-r)^2}\bigl(\,1-(rs)\lg\om_3',\om_4\rg\bigr),
\quad &(i,j)=(3,4),\\
\frac{1}{(s-r)^2}\bigl(\,\lg\om_4',\om_3\rg-(rs)\,\bigr),\quad
&(i,j)=(4,3).
\endcases\\
&=0, \quad (\text{\it since } \lg\om_3',\om_4\rg=r^{-1}s^{-1}, \
\& \ \lg\om_4',\om_3\rg=rs).
\endsplit
$$

Similarly, we can prove that the $(r,s)$-Serre relation (D5)  is preserved for ${\Cal
B'}^{\text{cop}}$.

\smallskip (I) and (II) together complete
the proof of the Proposition. \hfill\qed
\enddemo

\proclaim{Definition 2.4} \ For any two Hopf algebras
$\Cal A$ and $\Cal U$ paired by a skew-dual pairing $\lg\,,\rg$, one may
form the Drinfel'd (quantum) double $\Cal D(\Cal A,\Cal U)$ as in
\cite{12, 3.2} or \cite{13, 8.2}. This is a Hopf algebra whose
underlying vector space is $\Cal A\ot\Cal U$ with the tensor product
coalgebra structure and the algebra structure defined by
$$
(a\ot f)(a'\ot f')=\sum \lg S_{\Cal U}(f_{(1)}), a'_{(1)}\rg\lg
f_{(3)},a'_{(3)}\rg \,aa'_{(2)}\ot f_{(2)}f',\eqno{(8)}
$$
for $a, a'\in \Cal A$ and $f, f'\in\Cal U$, and whose antipode $S$
is given by
$$
S(a\ot f)=(1\ot S_{\Cal U}(f))(S_{\Cal A}(a)\ot 1). \eqno{(9)}
$$
\endproclaim

Clearly, both mappings $\Cal A\ni a\mapsto a\ot 1\in\Cal D(\Cal
A,\Cal U)$ and $\Cal U\ni f\mapsto 1\ot f\in\Cal D(\Cal A, \Cal
U)$ are injective Hopf algebra homomorphisms. Let us denote the
image $a\ot 1$ (resp. $1\ot f$) of $a$ (resp. $f$) in $\Cal D(\Cal
A,\Cal U)$ by $\hat a$ (resp. $\hat f$). By (8), we have the
following cross commutation relations between elements $\hat a$
(for $a\in\Cal A$) and $\hat f$ (for $f\in\Cal U$) in the algebra
$\Cal D(\Cal A,\Cal U)$:
$$\gather
\hat f\,\hat a=\sum\, \lg S_{\Cal U}(f_{(1)}), a_{(1)}\rg\,\lg
f_{(3)},a_{(3)}\rg\;\hat a_{(2)}\hat f_{(2)},\tag{10}\\
\sum\lg f_{(1)}, a_{(1)}\rg\,\hat f_{(2)}\,\hat a_{(2)}= \sum \hat
a_{(1)}\,\hat f_{(1)}\,\lg f_{(2)},a_{(2)}\rg.\tag{11}
\endgather
$$
In fact, as an algebra, the double $\Cal D(\Cal A,\Cal U)$ is the
universal algebra generated by the algebras $\Cal A$ and $\Cal U$
with cross relation (10) or equivalently (11).

\proclaim{Theorem 2.5} \ The two-parameter quantum group
$U=U_{r,s}(\frak g)$  is isomorphic to the Drinfel'd quantum
double $\Cal D(\Cal B, \Cal B')$, for $\frak g=\frak {sl}_n$,
$\frak {so}_{2n+1}$, $\frak {so}_{2n}$, or $\frak{sp}_{2n}$.
\endproclaim
\demo{Proof} \ Define a mapping $\varphi:\, \Cal D(\Cal B,\Cal
B')\lra U_{r,s}(\frak g)$ by
$$\gather
\varphi({\hat \om}_i^{\pm1})=\om_i^{\pm1},\; \quad \qquad
\varphi({{\hat \omega}_i}^{'\pm 1})={{\omega}'_i}^{\pm 1},\\
\varphi({\hat e}_i)=e_i, \qquad\quad\qquad \varphi({\hat f}_i)=f_i\,.
\endgather
$$

Note that by definition, $\varphi$ preserves the coalgebra
structures, the relations in $\Cal B$, and the relations in $\Cal
B'$.

We next verify that the cross relations in the double $\Cal
D(\Cal B,\Cal B')$ correspond to those in $U$.

By (11) and applying the comultiplication given in Proposition
1.1, we obtain the cross relations
$$\split
\lg{\om_j}^{'\pm1},\om_i^{\pm1}\rg\,{\hat\om_j}^{'\pm1}\,{\hat\om_i}^{\pm1}&=\hat\om_i^{\pm1}\,
{\hat\om_j}^{'\pm1}\lg
{\om_j}^{'\pm1},\om_i^{\pm1}\rg,\\
\lg \om_j',\om_i\rg\,\hat\om_j'\,\hat e_i&=\hat
e_i\,\hat\om_j'\,\lg \om_j',1\rg,\\
\lg 1,\om_j^{-1}\rg\,\hat
f_i\,\hat\om_j^{-1}&=\hat\om_j^{-1}\,\hat
f_i\,\lg\om_i',\om_j^{-1}\rg,\\
\lg f_j, e_i\rg\,\hat \om_j'+\lg 1,\om_i\rg\, \hat f_j\,\hat
e_i&=\hat e_i\,\hat f_j\,\lg \om_j',1\rg+\hat\om_i\,\lg f_j,
e_i\rg.
\endsplit
$$
That is,
$$\split
{\hat\om_j}^{'\pm1}\,{\hat\om_i}^{\pm1}&=\hat\om_i^{\pm1}\,{\hat\om_j}^{'\pm1},\\
\hat\om_j'\,\hat e_i\,{\hat\om_j}^{'-1}&=\hat e_i\,\lg
\om_j',\om_i\rg^{-1},\\
\hat\om_j\,\hat f_i\,\hat\om_j^{-1}&=\hat
f_i\,\lg\om_i',\om_j\rg^{-1},\\
[\,\hat e_i,\hat
f_j\,]&=\delta_{ij}\frac{\hat\om_i-\hat\om_i'}{r_i-s_i}.
\endsplit
$$
The coefficients $\lg \om_i',\om_j\rg^{-1}$ of $\hat
f_i$ in the third formula above coincide with those of the
formulae related to $f_i$ in $(X2)$ (here $X=A,\,B,\,C,\,D$,
resp.). By Lemma 1.2, we find that the coefficients $\lg
\om_j',\om_i\rg^{-1}$ of $\hat e_i$ in the second formula above do
coincide with those of the formulae related to $e_i$ in $(X3)$.
Hence, applying $\varphi$ gives the desired relations $(X2)$,
$(X3)$ and $(X4)$ in $U$.

As $U$ is generated by $e_i$, $f_i$, $\om_i^{\pm1}$ and
${\om_i'}^{\pm1}$ ($1\le i\le n$), the mapping $\varphi$ is
surjective. Since $U$ and $\Cal D(\Cal B,\Cal B')$ are
universal, subject to the same relations on essentially the
same generating set, $\varphi$ provides an isomorphism. \hfill\qed
\enddemo

\noindent {\bf Remarks.} \ (1) \ Up to now, we have completely solved
the compatibility problem on the defining relations of our
two-parameter quantum groups $U_{r,s}(\frak g)$ for $\frak
g=\frak{so}_{2n+1}$, $\frak{so}_{2n}$ and $\frak{sp}_{2n}$. This is done in two steps:
 the proof of Theorem 2.5 indicates that the cross
relations between $\Cal B$ and $\Cal B'$ are half of the relations
  (X1)---(X4), and the proof of Proposition 2.3
shows the remaining relations,  including the remaining half of
(X1)---(X4) and the $(r,s)$-Serre relations (X5)---(X7).

\medskip

 (2) \ Let $U^0=\Bbb
Q(r,s)[\om_1^{\pm1},\cdots,\om_n^{\pm1},{\om_1'}^{\pm1},\cdots,{\om_n'}^{\pm1}]$,
$U_0=\Bbb Q(r,s)[\om_1^{\pm1},\cdots,\om_n^{\pm1}]$, and
$U_0'=\Bbb Q(r,s)[{\om_1'}^{\pm1},\cdots,{\om_n'}^{\pm1}]$ denote
the Laurent polynomial subalgebras of $U_{r,s}(\frak
g)$, $\Cal B$,
 and $\Cal B'$ respectively. Clearly, $U^0=U_0U_0'=U_0'U_0$.
  Furthermore, let us
denote by $U_{r,s}(\frak n)$ $($resp. $U_{r,s}(\frak n^-)$\,$)$ the
subalgebra of $\Cal B$ $($resp. $\Cal B'$$)$ generated by $e_i$
$($resp. $f_i$$)$ for all $i\le n$. Thus, by definition, we have
$\Cal B=U_{r,s}(\frak n)\rtimes U_0$, and $\Cal B'=U_0'\ltimes
U_{r,s}(\frak n^-)$, so that the double $\Cal D(\Cal B,\Cal
B')\cong U_{r,s}(\frak n)\ot U^0\ot U_{r,s}(\frak n^-)$, as vector
spaces.

\medskip

(3) \ The above Theorem further implies the existence of the
following standard triangular decomposition of $U_{r,s}(\frak g)$,
which means that $U_{r,s}(\frak g)$ possesses highest weight
representation theory in the usual sense. In the corollary below,
we point out that the standard triangular decomposition structure
of $U_{r,s}(\frak g)$ in ``two-parameter" is a natural consequence
of its Drinfel'd double structure. This yields a direct approach
to the triangular decomposition structure of the ``one-parameter"
quantum groups $U_q(\frak g)$ of Drinfel'd-Jimbo type.  The two
original proofs by Lusztig \cite{16} and Rosso \cite{20} were
remarkable but nontrivial compared to ours.

\proclaim{Corollary 2.6} \ $U_{r,s}(\frak g)\cong U_{r,s}(\frak
n^-)\ot U^0\ot U_{r,s}(\frak n)$, as vector spaces. In particular,
it induces $U_q(\frak g)\cong U_q(\frak n^-)\otimes U_0\otimes
U_q(\frak n)$, as vector spaces.
\endproclaim
\demo{Proof} \ Define $\lg\,,\rg^-:\,\Cal B'\times\Cal B\lra\Bbb
Q(r,s)$ by $\lg b', b\rg^-:=\lg S(b'), b\rg$ for $b'\in\Cal B'$,
$b\in\Cal B$, which is the convolution inverse of the skew-dual
pairing $\lg\,,\rg$ in Proposition 2.3. By definition, it is
easily seen that its composition with the flip mapping $\tau$
yields a new skew-dual pairing
$\lg\,|\,\rg_0:=\lg\;,\,\rg^-\circ\tau: \Cal B\times\Cal B'\lra
\Bbb Q(r,s)$, given by $\lg b\,|\,b'\rg_0=\lg S(b'),
b\rg$.

Denote $\Cal D(\Cal B',\Cal B)$ the Drinfel'd double constructed
from the skew-dual pairing $\lg\,|\,\rg_0$. Then the following
mapping $\theta$ (see \cite{13, 8.2.1(26)}) establishes an
isomorphism between Hopf algebras $\Cal D(\Cal B,\Cal B')$ and
$\Cal D(\Cal B',\Cal B)$:
$$
\theta(b\ot b')=\sum \lg b_{(1)}',b_{(1)}\rg\, b_{(2)}'\ot
b_{(2)}\,\lg b_{(3)}|\,b'_{(3)}\rg_0, \qquad b\in\Cal B, \
b'\in\Cal B'.
$$
Hence, the composition $\theta\circ\varphi^{-1}$ yields the
required Hopf algebra isomorphism, where $\varphi$ is used in the
proof of Theorem 2.5.

If one takes $r=q$, $s=q^{-1}$, the above triangular decomposition
of $U_{r,s}(\frak g)$, modulo some identifications in $U^0$,
implies the ``one-parameter" case. \hfill\qed
\enddemo

Let $Q=\Bbb Z\Psi$ denote the root lattice and set
$Q^+=\sum_{i=1}^n\Bbb Z_{\ge0}\al_i$. Then for any
$\zeta=\sum_{i=1}^n\zeta_i\al_i\in Q$, we adopt the notation
$$
\om_\zeta=\om_1^{\zeta_1}\cdots\om_n^{\zeta_n}, \qquad
\om_\zeta'=(\om_1')^{\zeta_1}\cdots(\om_n')^{\zeta_n}.
$$
We obtain a $Q$-graded structure on $U_{r,s}(\frak g)$ as another Corollary
of Proposition 2.3 and Theorem 2.5. This will be useful
for the representation theory discussed later.

\proclaim{Corollary 2.7} \ For any
$\zeta=\sum_{i=1}^n\zeta_i\al_i\in Q$, the defining relations
$(X2)$ and $(X3)$ in $U_{r,s}(\frak g)$ take the form:
$$\gather
\om_{\zeta}\,e_i\,\om_{\zeta}^{-1}=\lg \om_i',\om_\zeta\rg\,e_i,
\qquad
\om_{\zeta}\,f_i\,\om_{\zeta}^{-1}=\lg \om_i',\om_\zeta\rg^{-1}f_i,\\
{\om_{\zeta}'}\,e_i\,{\om_{\zeta}'}^{-1}=\lg \om_\zeta',
\om_i\rg^{-1} e_i,\qquad\quad
\om_{\zeta}'\,f_i\,{\om_{\zeta}'}^{-1}=\lg
\om_\zeta',\om_i\rg\,f_i.
\endgather
$$
$U_{r,s}(\frak n^\pm)=\bigoplus_{\eta\in
Q^+}U_{r,s}^{\pm\eta}(\frak n^\pm)$ is then $Q^\pm$-graded with
$$
U_{r,s}^{\eta}(\frak n^\pm)=\left\{\,a\in U_{r,s}(\frak
n^\pm)\;\left|\; \om_\zeta\,a\,\om_\zeta^{-1}=\lg
\om_\eta',\om_\zeta\rg\,a, \ \om_\zeta'\,a\,{\om_\zeta'}^{-1}=\lg
\om_\zeta',\om_\eta\rg^{-1} \,a\,\right\}\right.,\eqno{(12)}
$$
for $\eta\in Q^+\cup Q^-$.

\medskip
Furthermore, $U=\bigoplus_{\eta\in Q}U_{r,s}^\eta(\frak g)$ is
$Q$-graded with
$$
\split U_{r,s}^\eta(\frak g)&=\left\{\,\left.\sum
F_\al\om_{\mu}'\om_\nu E_\be\in U\; \right|\;
\om_\zeta\,(F_\al\om_{\mu}'\om_\nu E_\be)\,\om_{\zeta}^{-1}=
\lg \om'_{\beta-\al},\om_\zeta\rg\,F_\al\om_\mu'\om_\nu E_\be,\right.\\
&\ \left.\om_\zeta'\,(F_\al\om_\mu'\om_\nu
E_\be)\,{\om_{\zeta}'}^{-1}= \lg
\om_\zeta',\om_{\beta-\al}\rg^{-1}\,F_\al\om_\mu'\om_\nu E_\be,\;
\text{\it with } \; \beta-\al=\eta\right\},
\endsplit\tag{13}
$$
where $F_\al$ $($resp. $E_\be$$)$ runs over monomials
$f_{i_1}{\cdots} f_{i_l}$ $($resp. $e_{j_1}{\cdots} e_{j_m}$$)$
such that $\al_{i_1}+{\cdots}+\al_{i_l}=\al$ $($resp.
$\al_{j_1}+{\cdots}+\al_{j_m}=\be$$)$.\hfill\qed
\endproclaim

Let $\lg\,|\,\rg_0:\,\Cal B\times \Cal B'\lra \Bbb Q(r,s)$ denote
the skew-dual pairing introduced in the proof of Corollary 2.6.
Then we have \proclaim{Definition 2.8} \ The bilinear form
$\lg\,,\,\rg_U$ on $U_{r,s}(\frak g)\times U_{r,s}(\frak g)$
defined by
$$\split
\lg F_\al\om_{\mu}'\om_\nu E_\be,F_\theta\om_{\sigma}'\om_\delta
E_\gamma\rg_U&=\lg S(\om_\nu
E_\beta)\,|\,F_\theta\om_\sigma'\rg_0\lg\om_\delta
E_\gamma\,|\,S(F_\al\om_\mu')\rg_0\\
&=\lg\om_\nu\,|\,\om_\sigma'\rg_0^{-1}\lg\om_\delta\,|\,\om_\mu'\rg_0^{-1}\lg
S(E_\be)\,|\,F_\theta\rg_0\lg E_\gamma\,|\,S(F_\al)\rg_0\\
&=\lg\om_\sigma',\om_\nu\rg\lg\om_\mu',\om_\delta\rg\lg
F_\theta,E_\be\rg\lg S^2(F_\al),E_\gamma\rg
\endsplit\tag{14}
$$
is called the Rosso form of the two-parameter quantum group
$U_{r,s}(\frak g)$.
\endproclaim

\proclaim{Proposition 2.9} \ The Rosso form $\lg\,,\,\rg_U$ on
$U_{r,s}(\frak g)\times U_{r,s}(\frak g)$ is $\text{\rm
ad}_l$-invariant, that is, $$ \lg\text{\rm
ad}_l(a)\,b,\,c\rg_U=\lg b,\,\text{\rm ad}_l(S(a))\,c\rg_U, \qquad
a,\,b,\,c\in U_{r,s}(\frak g).
$$
\endproclaim
\demo{Proof} \ By Corollary 2.6, $U_{r,s}(\frak g)\cong\Cal D(\Cal
B',\Cal B)$ with respect to $\lg\,|\,\rg_0$. By \cite{13,
Proposition 8.12}, the Rosso form $\lg\cdot,\cdot\rg_U$ on the
quantum double $\Cal D(\Cal B',\Cal B)$ is
$\text{ad}_l$-invariant. \hfill\qed
\enddemo

\bigskip\bigskip
\heading{$\bold 3$. \ Lusztig's Symmetry } \endheading
\bigskip

We define in this section the Lusztig's symmetries for the
two-parameter quantum groups $U_{r,s}(\frak g)$ we have defined in
Section 1. A remarkable feature of these symmetries in the
``two-parameter" cases is their existence between these quantum
groups and the so-called associated objects only as $\Bbb
Q$-isomorphisms rather than as $\Bbb Q(r,s)$-automorphisms as
usual in the ``one-parameter" cases. When $\text{\rm rank}\,(\frak
g)>2$, we find that we can construct the Lusztig's symmetries of
$U_{r,s}(\frak g)$ into its associated quantum group
$U_{s^{-1},r^{-1}}(\frak g)$ (as defined below)
 if and only if $rs=1$. In this case, if we set $r=q$
and $s=q^{-1}$, the Lusztig's symmetries turn out to be $\Bbb
Q(q)$-automorphisms of $U_{q,q^{-1}}(\frak g)$, which in
particular induce the usual Lusztig's symmetries defined on
quantum groups $U_q(\frak g)$ of Drinfel'd-Jimbo type.

\proclaim{Theorem 3.1} \ $\text{\rm(i)}$ \ When $\text{\rm
rank}\,(\frak g)=2$, for $\frak g=\frak {sl}_3$, $\frak {sp}_4$,
or $\frak {so}_5$, the Lusztig's symmetries exist between
$U_{r,s}(\frak g)$ and its associated quantum group
$U_{s^{-1},r^{-1}}(\frak g)$.

$\text{\rm(ii)}$ \ When $\text{\rm rank}\,(\frak g)>2$ for any
type of $\frak g$, the two-parameter quantum group $U_{r,s}(\frak
g)$ has the Lusztig's symmetries if and only if it is of the form
$U_{q,q^{-1}}(\frak g)$, where $r=q$, $s=q^{-1}$. In particular,
in this case, each Lusztig's symmetry induces a usual Lusztig's
symmetry on $U_q(\frak g)$.
\endproclaim

Before giving the proof, we need to make some general preliminary
remarks. We first observe that the pairing $\lg\,,\rg$ in
Proposition 2.3 plays a role in locating the structure constants
of $U_{r,s}(\frak g)$. When it is necessary to emphasize this
point, we denote the two-parameter quantum groups by
$(U_{r,s}(\frak g), \lg\,,\rg)$. Now we call
$(U_{s^{-1},r^{-1}}(\frak g), \lg\,|\,\rg)$ the associated quantum
group corresponding to $(U_{r,s}(\frak g), \lg\,,\rg)$, where the
pairing $\lg\om_i'|\,\om_j\rg$ is defined via substituting $(r,
s)$ by $(s^{-1}, r^{-1})$ in the defining formula for $\lg
\om_i',\om_j\rg$. Lemma 1.2 guarantees that
$$
\lg \om_i'|\,\om_j\rg=\lg\om_j',\om_i\rg.\eqno(1)
$$

In order to define the Lusztig's symmetries, we introduce the
notion of divided-power elements: for any nonnegative integer
$k\in \Bbb N$, set
$$\gather
[k]=\frac{r^k-s^k}{r-s},\qquad [k]!=[1][2]\cdots [k],\\
\lg k\rg=\frac{s^{-k}-r^{-k}}{s^{-1}-r^{-1}},\qquad \lg k\rg!=\lg
1\rg\lg 2\rg\cdots\lg k\rg,
\endgather
$$
and for any element $x\in U_{r,s}(\frak g)$ (or
$U_{s^{-1},r^{-1}}(\frak g)$), we define two kinds of divided-power
elements:
$$
x^{[k]}=x^{k}/[k]!\,,\qquad x^{(k)}=x^{k}/\lg k\rg!\,.
$$

\proclaim{Definition 3.2} \ To every $i$, $i=1,\cdots,n$, there
corresponds a $\Bbb Q$-linear mapping $\Cal T_i:\, (U_{r,s}(\frak
g), \lg\,,\rg)\lra (U_{s^{-1},r^{-1}}(\frak g), \lg\,| \,\rg)$
defined on generators $\om_j$, $\om_j'$, $e_j$, $f_j$ as follow:
$$\gather
\Cal T_i(\om_j)=\om_j\,\om_i^{-a_{ij}},\qquad \Cal
T_i(\om_j')=\om_j'\,{\om_i'}^{-a_{ij}},\\
\Cal T_i(e_i)=-\,{\om_i'}^{-1}f_i,\qquad \Cal
T_i(f_i)=-(r_is_i)\,e_i\,\om_i^{-1};\\
\endgather
$$
when $\frak g$ is of type $A$, $C$, or $D$,
$$\gather \Cal
T_i(e_j)=\sum_{\nu=0}^{-a_{ij}}(-1)^{\nu}(rs)^{\frac{\nu}2(-a_{ij}-\nu)}\lg\om_j',\om_i\rg^{-\nu}
\lg\om_i',\om_i\rg^{\frac{\nu}{2}(1{+}a_{ij})}
e_i^{(\nu)}e_j\,e_i^{(-a_{ij}-\nu)},\; i\ne j,\\
\Cal
T_i(f_j)=(r_js_j)^{\delta_{ij}^+}\sum_{\nu=0}^{-a_{ij}}(-1)^{\nu}(rs)^{\frac{\nu}2(-a_{ij}-\nu)}
\lg\om_i',\om_j\rg^\nu\lg\om_i',\om_i\rg^{-\frac{\nu}{2}(1{+}a_{ij})}
f_i^{(-a_{ij}-\nu)}f_j\,f_i^{(\nu)},\; i\ne j;\\
\endgather
$$
and when $\frak g$ is of type $B$,
$$\gather \Cal
T_i(e_j)=\sum_{\nu=0}^{-a_{ij}}(-1)^{\nu}(rs)^{\frac{\nu}2(-a_{ij}-\nu)}\lg\om_j',\om_i\rg^{\nu}
\lg\om_i',\om_i\rg^{-\frac{\nu}{2}(1{+}a_{ij})}
e_i^{(-a_{ij}-\nu)}e_j\,e_i^{(\nu)},\; i\ne j,\\
\Cal
T_i(f_j)=(r_js_j)^{\delta_{ij}^-}\sum_{\nu=0}^{-a_{ij}}(-1)^{\nu}(rs)^{\frac{\nu}2(-a_{ij}-\nu)}
\lg\om_i',\om_j\rg^{-\nu}\lg\om_i',\om_i\rg^{\frac{\nu}{2}(1{+}a_{ij})}
f_i^{(\nu)}f_j\,f_i^{(-a_{ij}-\nu)},\; i\ne j.\\
\endgather
$$
Furthermore  $\Cal T_i(r)=s^{-1}$ and $\Cal T_i(s)=r^{-1}$ so that
$\Cal T_i(x^{[k]})=(\Cal T_i(x))^{(k)}$.  Here $(a_{ij})$ is the
Cartan matrix of the classical simple Lie algebra $\frak g$, and
for any $i\ne j$,
$$\gather
\delta_{ij}^+=\cases 2, & \text{ if }\ i<j,\ \&\ a_{ij}\ne0, \\
1, & \text{ otherwise }.\\
\endcases
\qquad\quad \delta_{ij}^-=\cases 2, & \text{ if }\quad i>j,\ \&\ a_{ij}\ne0, \\
1, & \text{ otherwise }.\\
\endcases\\
 e_i^{(k)}=e_i^{k}/\lg k\rg_i!, \qquad f_i^{(k)}=f_i^{k}/\lg
k\rg_i!,\qquad \lg
k\rg_i=\frac{s_i^{-k}-r_i^{-k}}{s_i^{-1}-r_i^{-1}}.
\endgather
$$
\endproclaim

\noindent {\bf Remarks.} \ (1) \ The defining formulas of
the actions of $\Cal T_i$ on $e_j$ and $f_j$ when $i\ne j$
(in the cases of types $A$, $C$, and $D$) can be interpreted
as the right- and the left-adjoint actions of the
divided-power operators $(\text{ad}_re_i)^{(-a_{ij})}$ and
$(\text{ad}_lf_i)^{(-a_{ij})}$ in the co-opposite Hopf algebra
$U_{r,s}(\frak g)^{\text{cop}}$:
$$
\Cal T_i(e_j)=(\text{ad}_re_i)^{(-a_{ij})}(e_j), \qquad \Cal
T_i(f_j)=(r_js_j)^{\delta_{ij}^+}(\text{ad}_lf_i)^{(-a_{ij})}(f_j).
$$
This phenomenon for the ``two-parameter" cases is interesting as
it implies that $U_{r,s}(\frak g)^{\text{cop}}$ $\cong
U_{s^{-1},r^{-1}}(\frak g)$ as Hopf algebras. While for the case
of type $B$, it can be roughly viewed as the dual case of type
$C$.

According to Definition 3.2, we can also consider the Lusztig's
symmetries $\Cal T_i$ defined for $U_{r,s}(\frak g)$ as $\Bbb
Q$-automorphisms of $U_{r,s}(\frak g)$ into itself. However, the
images of $\Cal T_i$ should be contained ``locally" (here we mean
only in each rank $2$ size) in its associated quantum group
$U_{s^{-1},r^{-1}}(\frak g)$. This fact will be clear, in light of
the proofs of Lemmas 3.3 --- 3.5 below.

\smallskip

(2) \ When $r=s^{-1}=q$, the Lusztig's symmetries $\Cal T_i$ are
$\Bbb Q(q)$-automorphisms of $U_{q,q^{-1}}(\frak g)$. This is
clear if we extend the action of $\Cal T_i$ on $U_{r,s}(\frak g)$
algebraically, via an identification, we then get the usual
Lusztig's symmetries defined on the quantum groups $U_q(\frak g)$
of Drinfel'd-Jimbo type.

\medskip In the following, we  consider the special cases in
rank $2$: $A_2$, $B_2$, or $C_2$.

\proclaim {Lemma 3.3} \ Assume that $\frak g=\frak {sl}_3$,
$\frak{so}_5$, or $\frak {sp}_4$. Then $\Cal T_i$ $($$i=1,2$$)$
preserve the defining relations $(X1)$---$(X3)$ of $(U_{r,s}(\frak
g), \langle\,,\rangle)$ into its associated object $(U_{s^{-1},
r^{-1}}(\frak g), \langle\,|\,\rangle)$, where $X=A$, $B$, or $C$.
\endproclaim
\demo{Proof} \ In case $(A_2)$: we have
$$\gather
\lg\om_1',\om_1\rg=rs^{-1}=\lg\om_1'|\,\om_1\rg, \qquad
\lg\om_1',\om_2\rg=r^{-1}=\lg\om_2'|\,\om_1\rg,\\
\lg\om_2',\om_1\rg=s=\lg\om_1'|\,\om_2\rg,\qquad
\lg\om_2',\om_2\rg=rs^{-1}=\lg\om_2'|\,\om_2\rg.
\endgather
$$
In case $(B_2)$: we have
$$\gather
\lg\om_1',\om_1\rg=(rs^{-1})^2=\lg\om_1'|\,\om_1\rg,\qquad
\lg\om_1',\om_2\rg=r^{-2}=\lg\om_2'|\,\om_1\rg,\\
\lg\om_2',\om_1\rg=s^2=\lg\om_1'|\,\om_2\rg,\qquad
\lg\om_2',\om_2\rg=rs^{-1}=\lg\om_2'|\,\om_2\rg.
\endgather
$$
In case $(C_2)$: we have
$$\gather
\lg\om_1',\om_1\rg=rs^{-1}=\lg\om_1'|\,\om_1\rg,\qquad
\lg\om_1',\om_2\rg=r^{-2}=\lg\om_2'|\,\om_1\rg,\\
\lg\om_2',\om_1\rg=s^2=\lg\om_1'|\,\om_2\rg,\qquad
\lg\om_2',\om_2\rg=r^2s^{-2}=\lg\om_2'|\,\om_2\rg.
\endgather
$$

We need to show that $\Cal T_1,\,\Cal T_2$ preserve the defining
relations $(X1)$---$(X3)$. $(X1)$ are automatically satisfied.

To check $(X2)$ \& $(X3)$ we first remark that in the rank 2 cases, we have $\Cal T_k(\lg
\om_i',\om_j\rg)=\lg \Cal T_k(\om_i'),\Cal T_k(\om_j)\rg=\lg
\om_j',\om_i\rg=\lg \om_i'|\,\om_j\rg$, for $i,\,j,\,k\in
\{1,\,2\}$. This fact ensures that $\Cal T_k$ ($k=1,\,2$) preserve
$(X2)$ and $(X3)$, that is
$$
\gather \Cal T_k(\om_j)\Cal T_k(e_i)\Cal
T_k(\om_j)^{-1}=\lg\om_i'|\,\om_j\rg\,\Cal T_k(e_i), \qquad \Cal
T_k(\om_j)\Cal T_k(f_i)\Cal
T_k(\om_j)^{-1}=\lg\om_i'|\,\om_j\rg^{-1}\Cal T_k(f_i),
\\
\Cal T_k(\om_j')\Cal T_k(e_i)\Cal
T_k(\om_j')^{-1}=\lg\om_j'|\,\om_i\rg^{-1}\Cal T_k(e_i),\qquad
\Cal T_k(\om_j')\Cal T_k(f_i)\Cal
T_k(\om_j')^{-1}=\lg\om_j'|\,\om_i\rg\,\Cal T_k(f_i).\\
\endgather
$$
 All identities follow from the first one. \hfill\qed
\enddemo

\proclaim{Lemma 3.4} \ With the same assumption as in Lemma 3.3,
$\Cal T_i$ $($$i=1,\,2$$)$ preserves the defining relations $(X4)$
of $(U_{r,s}(\frak g), \langle\,,\rangle)$ into its associated
object $(U_{s^{-1}, r^{-1}}(\frak g), \langle\,|\,\rangle)$, for
$X=A,\,B,\,C$.
\endproclaim
\demo{Proof} \ For $i=1,\,2$, we have
$$\split
[\,\Cal T_i(e_i),\Cal T_i(f_i)\,]
&=(r_is_i)\,{\om_i'}^{-1}\bigl(f_i\,e_i-\lg\om_i',\om_i\rg^{-1}\lg\om_i',\om_i\rg\,e_i\,f_i\bigr)\,\om_i^{-1}\\
&=\frac{\om_i^{-1}-{\om_i'}^{-1}}{s_i^{-1}-r_i^{-1}}=\frac{\Cal
T_i(\om_i)-\Cal T_i(\om_i')}{s_i^{-1}-r_i^{-1}}\\
&=\Cal T_i([\,e_i, f_i\,])\in U_{s^{-1},r^{-1}}(\frak g).
\endsplit
$$

When $a_{ij}=-1$, we have
$$
\split [\,\Cal T_i(e_j),\Cal T_i(f_j)\,]&=
(r_js_j)^{\delta_{ij}^+}\,[\,e_j\,e_i-\lg\om_j',\om_i\rg^{-1}e_i\,e_j,
f_i\,f_j-\lg\om_i',\om_j\rg\,f_j\,f_i\,]\\
&=(r_js_j)^{\delta_{ij}^+}\left\{f_i[e_j,f_j]e_i+e_j[e_i,f_i]f_j-\lg\om_j',
\om_i\rg^{-1}\,([e_i,f_i]f_je_j+e_if_i[e_j,f_j])\right.\\
&\quad -\lg \om_i',\om_j\rg([e_j,f_j]f_ie_i+e_jf_j[e_i,f_i])\\
&\quad \left.+\lg\om_j',\om_i\rg^{-1}\lg \om_i',\om_j\rg\,(e_i[e_j,f_j]f_i+f_j[e_i,f_i]e_j)\right\}\\
&=\frac{(r_js_j)^{\delta_{ij}^+}}{r_j-s_j}\frac{1}{r_i-s_i}\,\left\{(\om_j-\om_j')\,\bigl[\,
\om_i\lg\om_j',\om_i\rg^{-1}-\om_i'\lg\om_i',\om_j\rg\right.\\
&\quad-(\lg\om_j',\om_i\rg^{-1}+\lg\om_i',\om_j\rg)(\om_i-\om_i')\,\bigr]\\
&\quad
\left.+\lg\om_j',\om_i\rg^{-1}\lg\om_i',\om_j\rg\bigl(\om_j\lg\om_i',\om_j\rg^{-1}
-\om_j'\lg\om_j',\om_i\rg\bigr)(\om_i-\om_i')\right\}\\
&=\frac{(r_js_j)^{\delta_{ij}^+}\bigl(\lg\om_j',\om_i\rg^{-1}-\lg\om_i',\om_j\rg\bigr)}{(r_j-s_j)(r_i-s_i)}
(\om_i\om_j-\om_i'\om_j')\\
&=\cases \frac{\Cal T_i(\om_j)-\Cal T_i(\om_j')}{s^{-1}-r^{-1}}, &
\text{\it for } \ (i, j)=(1,
2), \; i.e., \; \text{\it in type } \; A_2\\
\frac{\Cal T_i(\om_j)-\Cal T_i(\om_j')}{s^{-1}-r^{-1}}, &
\text{\it for } \ (i, j)=(2, 1), \; i.e., \; \text{\it in types }
\; A_2, \; C_2
\endcases\\
&=\Cal T_i([e_j, f_j]).
\endsplit
$$

When $a_{12}=-2$ in type $C_2$, we need to prove
$$
[\,\Cal T_1(e_2),\Cal T_1(f_2)\,]=
\frac{\om_2\om_1^2-\om_2'{\om_1'}^2}{s^{-2}-r^{-2}} =\frac{\Cal
T_1(\om_2)-\Cal T_1(\om_2')}{s_2^{-1}-r_2^{-1}}=\Cal T_1([e_2,
f_2])\in U_{s^{-1},r^{-1}}(\frak g).$$ To this end, let us denote
$$\gather
E_{12}=e_{\al_1+\al_2}=\text{ad}_le_1(e_2)=e_1e_2-s^2e_2e_1, \\
F_{12}=f_{\al_1+\al_2}=\text{ad}_rf_1(f_2)=f_2f_1-r^2f_1f_2,\\
E_{112}=e_{2\al_1+\al_2}=\text{ad}_le_1^2(e_2)=e_1E_{12}-rsE_{12}e_1,
\\
F_{112}=f_{2\al_1+\al_2}=\text{ad}_rf_1^2(f_2)=F_{12}f_1-rsf_1F_{12}.
\endgather
$$
A direct calculation shows
$$\gather
[e_1, F_{12}]=-(r+s)\om_1f_2, \qquad [e_2,
F_{12}]=f_1\om_2',\\
[E_{12}, f_1]=-(r+s)e_2\om_1', \qquad [E_{12}, f_2]=\om_2e_1,\\
[E_{12}, F_{12}]=\frac{\om_1\om_2-\om'_1\om'_2}{r-s}.
\endgather
$$
Using the Leibniz rule and the above results, we can get
$$\split
[E_{112}, F_{112}]&=[e_1E_{12}-rsE_{12}e_1,
F_{12}f_1-rsf_1F_{12}]\\
&=[e_1, F_{12}]f_1E_{12}+e_1[E_{12}, F_{12}]f_1+F_{12}[e_1,
f_1]E_{12}+e_1F_{12}[E_{12}, f_1]\\
&\quad -(rs)([E_{12}, F_{12}]f_1e_1+E_{12}[e_1,
F_{12}]f_1+F_{12}[E_{12}, f_1]e_1+E_{12}F_{12}[e_1, f_1])\\
&\quad-(rs)([e_1, f_1]F_{12}E_{12}+f_1[e_1,
F_{12}]E_{12}+e_1[E_{12}, f_1]F_{12}+e_1f_1[E_{12},F_{12}])\\
&\quad+(rs)^2(E_{12}[e_1, f_1]F_{12}+E_{12}f_1[e_1, F_{12}]+
[E_{12}, f_1]F_{12}e_1+f_1[E_{12}, F_{12}]e_1)\\
&=\frac{(r+s)^2}{r^2-s^2}(\om_1^2\om_2-\om_1^{'2}\om'_2).
\endsplit
$$

Observe that $\Cal T_1(e_2)=\frac1{\lg 2\rg_1}E_{112}$, $\Cal
T_1(f_2)=\frac1{\lg 2\rg_1}F_{112}$, where $\lg
2\rg_1=(r{+}s)(rs)^{-1}$. We then arrive at the required equality
above.

A similar argument is used to prove the result for type $B_2$. \hfill\qed
\enddemo

\proclaim{Lemma 3.5} \ With the same assumption as in Lemma 3.3,
$\Cal T_i$ $($$i=1,\,2$$)$ preserves the $(r,s)$-Serre relations
$(X5)$, $(X6)$ and $(X7)$ of $(U_{r,s}(\frak g),
\langle\,,\rangle)$ into its associated object $(U_{s^{-1},
r^{-1}}(\frak g), \langle\,|\,\rangle)$, for $X=A$, $B$, $C$.
\endproclaim
\demo{Proof} \ We do not need to consider the relation $(X5)$ as it doesn't appear for the rank
$2$ cases.

For $A_2$: \ Consider the degree $2$ $(r,s)$-Serre relation
$$e_1^2e_2-(r+s)e_1e_2e_1+(rs)e_2e_1^2=0.\eqno(2)$$
Note that
$$\gather
\Cal T_1(e_2)\Cal T_1(e_1)=s\Cal T_1(e_1)\Cal
T_1(e_2)-r^{-1}e_2,\qquad e_2\Cal T_1(e_1)=r\Cal T_1(e_1)e_2;\\
\Cal T_2(e_1)\Cal T_2(e_2)=r^{-1}\Cal T_2(e_2)\Cal
T_2(e_1)-r^{-1}e_1,\qquad e_1\Cal T_2(e_1)=s\Cal T_2(e_1)e_1.
\endgather
$$
Acting with $\Cal T_i$ ($i=1,\,2$) algebraically on the left-hand side
of (2), we can get
$$\gather
\Cal T_1(e_1)^2\Cal T_1(e_2)-(r^{-1}+s^{-1})\Cal T_1(e_1)\Cal
T_1(e_2)\Cal T_1(e_1)+(r^{-1}s^{-1})\Cal T_1(e_2)\Cal
T_1(e_1)^2=0,\\
\Cal T_2(e_1)^2\Cal T_2(e_2)-(r^{-1}+s^{-1})\Cal T_2(e_1)\Cal
T_2(e_2)\Cal T_2(e_1)+(r^{-1}s^{-1})\Cal T_2(e_2)\Cal
T_2(e_1)^2=0,
\endgather
$$
that is, $\Cal T_i$ ($i=1,\,2$) preserves (2) for $A_2$.

Now consider another degree $2$ $(r,s)$-Serre relation
$$e_2^2e_1-(r_2^{-1}+s_2^{-1})e_2e_1e_2+(r_2^{-1}s_2^{-1})e_1e_2^2=0,\eqno(3)$$
which holds for both $A_2$ and $C_2$. Note that
$$
\Cal T_2(e_1)\Cal T_2(e_2)=r_2^{-1}\Cal T_2(e_2)\Cal
T_2(e_1)-r_2^{-1}e_1, \qquad \Cal T_2(e_2)e_1=s_2e_1\Cal T_2(e_2).
$$
Acting with $\Cal T_2$ algebraically on the left-hand side of (2), we
can easily get
$$
\Cal T_2(e_2)^2\Cal T_2(e_1)-(r_2+s_2)\Cal T_2(e_2)\Cal
T_2(e_1)\Cal T_2(e_2)+(r_2s_2)\Cal T_2(e_1)\Cal T_2(e_2)^2=0,
$$
that is, $\Cal T_2$ preserves (3) in both cases $A_2$ and $C_2$.

In the case of $A_2$, let $\Cal T_1$ act algebraically on the
left-hand side of (3). Because
$$\Cal T_1(e_2)\Cal T_1(e_1)=s\Cal
T_1(e_1)\Cal T_1(e_2)-r^{-1}e_2,\qquad e_2\Cal T_1(e_2)=r^{-1}\Cal
T_1(e_2)e_2,$$ we can easily get that
$$
\Cal T_1(e_2)^2\Cal T_1(e_1)-(r+s)\Cal T_1(e_2)\Cal T_1(e_1)\Cal
T_1(e_2)+(rs)\Cal T_1(e_1)\Cal T_1(e_2)^2=0.
$$
That is, $\Cal T_1$ preserves (3) for $A_2$, as well.

For $C_2$: \ Note that $[e_1^{(2)},
f_1]=\frac{s\om_1-r\om_1'}{r-s}e_1$. So by direct calculation, we
get
$$\gather
\Cal T_1(e_2)\Cal T_1(e_1)=s_2\Cal T_1(e_1)\Cal
T_1(e_2)+(rs)^{-1}\Cal T_1'(e_2),\tag{4}\\
\Cal T_1(e_2)\Cal T_1'(e_2)=r^2\Cal T_1'(e_2)\Cal T_1(e_2),\quad
(\text{by Lemma 3.6})
\endgather
$$
where $\Cal
T_1(e_2)=e_2e_1^{(2)}-s^{-1}e_1e_2e_1+r^{-1}s^{-3}e_1^{(2)}e_2$,
$\Cal T_1'(e_2)=(\text{ad}_le_1)(e_2)=e_1e_2-s^2e_2e_1$.

Let $\Cal T_1$ act algebraically on the left-hand side of (3), we
get
$$
\Cal T_1(e_2)^2\Cal T_1(e_1)-(r_2+s_2)\Cal T_1(e_2)\Cal
T_1(e_1)\Cal T_1(e_2)+(r_2s_2)\Cal T_1(e_1)\Cal T_1(e_2)^2=0.
$$
That is, $\Cal T_1$ preserves (3) for $C_2$.

Now consider the degree $3$ $(r,s)$-Serre relation
$$
e_1^3e_2-(r^2+rs+s^2)\,e_1^2e_2e_1+(rs)(r^2+rs+s^2)\,e_1e_2e_1^2-(rs)^3e_2e_1^3=0.\eqno(5)
$$
Observing that
$$\gather
e_2\Cal T_1(e_1)=r^2\Cal T_1(e_1)e_2,\\
\Cal T_1'(e_2)\Cal T_1(e_1)=(rs)\Cal T_1(e_1)\Cal
T_1'(e_2)+r^{-1}s(r+s)\,e_2,
\endgather
$$
and using (4), it is easy to check that $\Cal T_1$ preserves (5).

On the other hand, note that
$$\gather
\Cal T_2(e_1)\Cal T_2(e_2)=r^{-2}\Cal T_2(e_2)\Cal T_2(e_1)-r^{-2}e_1,\\
e_1\Cal T_2(e_1)^2-s(r+s)\Cal T_2(e_1)e_1\Cal T_2(e_1)+rs^3\Cal
T_2(e_1)^2e_1=0,\quad (\text{by Lemma 3.6})
\endgather
$$
which ensure $\Cal T_2$ preserves (5).

Dually, we can verify $\Cal T_i$ ($i=1,\,2$) preserves the
$(r,s)$-Serre relations $(X7)$ for $X=A_2,\,C_2$.

We can prove the result for type $B_2$ in a similar way.
\hfill\qed
\enddemo

\proclaim{Lemma 3.6} \ For type $C_2$, the following
identities hold.
$$\gather
\Cal T_1(e_2)\Cal T_1'(e_2)=r^2\Cal T_1'(e_2)\Cal
T_1(e_2),\tag{6}\\
e_1\Cal T_2(e_1)^2-s(r+s)\Cal T_2(e_1)e_1\Cal T_2(e_1)+rs^3\Cal
T_2(e_1)^2e_1=0,\tag{7}
\endgather
$$
where $\Cal
T_1(e_2)=e_2e_1^{(2)}-s^{-1}e_1e_2e_1+r^{-1}s^{-3}e_1^{(2)}e_2$,
$\Cal T_1'(e_2)=(\text{ad}_le_1)(e_2)=e_1e_2-s^2e_2e_1$, and $\Cal
T_2(e_1)=e_1e_2-r^2e_2e_1$.
\endproclaim
\demo{Proof} \ Writing
$\Delta=r^2+rs+s^2$, and using (C6), we get
$$\split
(rs^3&e_2e_1^2-s(r+s)e_1e_2e_1+e_1^2e_2)(e_1e_2-s^2e_2e_1)\\
&\quad-
r^2(e_1e_2-s^2e_2e_1)(rs^3e_2e_1^2-s(r+s)e_1e_2e_1+e_1^2e_2)\\
\endsplit
$$
$$\split
&=rs^3e_2e_1^3e_2-rs^5e_2e_1^2e_2e_1-s(r+s)e_1e_2e_1^2e_2\\
&\quad+s^3(r+s)e_1e_2e_1e_2e_1+e_1^2e_2e_1e_2-s^2e_1^2e_2^2e_1\\
&\quad-\bigl[(rs)^3e_1e_2^2e_1^2-r^2s(r+s)e_1e_2e_1e_2e_1+r^2e_1e_2e_1^2e_2\\
&\quad-r^3s^5e_2e_1e_2e_1^2+r^2s^3(r+s)e_2e_1^2e_2e_1-(rs)^2e_2e_1^3e_2\bigr]\\
&=rs^2(r{+}s)e_2e_1^3e_2-rs^3\Delta e_2(e_1^2e_2e_1)-\Delta
(e_1e_2e_1^2)e_2+s(r{+}s)(r^2{+}s^2)e_1e_2e_1e_2e_1\\
&\quad+e_1^2e_2e_1e_2-s^2e_1^2e_2^2e_1+r^3s^5e_2e_1e_2e_1^2-(rs)^3e_1e_2^2e_1^2\\
&=rs^2(r{+}s)e_2e_1^3e_2-rs^3e_2e_1^3e_2-r^2s^4\Delta
e_2e_1e_2e_1^2+r^4s^6e_2^2e_1^3\\
&\quad-(rs)^{-1}\Delta
e_1^2e_2e_1e_2+(rs)^{-1}e_1^3e_2^2-(rs)^2e_2e_1^3e_2+s(r{+}s)(r^2{+}s^2)e_1e_2e_1e_2e_1\\
&\quad+e_1^2e_2e_1e_2-s^2e_1^2e_2^2e_1+r^3s^5e_2e_1e_2e_1^2-(rs)^3e_1e_2^2e_1^2\\
&=-r^2s^4(r^2{+}s^2)(e_2e_1e_2)e_1^2+r^4s^6e_2^2e_1^3-(rs)^{-1}(r^2{+}s^2)e_1^2(e_2e_1e_2)\\
&\quad+(rs)^{-1}e_1^3e_2^2+s(r{+}s)(r^2{+}s^2)e_1(e_2e_1e_2)e_1-s^2e_1^2e_2^2e_1-(rs)^3e_1e_2^2e_1^2\\
&=-r^2s^4\bigl[(rs)^2e_2^2e_1^3+e_1e_2^2e_1^2\bigr]+
r^4s^6e_2^2e_1^3-(rs)^{-1}\bigl[(rs)^2e_1^2e_2^2e_1+e_1^3e_2^2\bigr]\\
&\quad+(rs)^{-1}e_1^3e_2^2+s(r+s)\bigl[(rs)^2e_1e_2^2e_1^2+e_1^2e_2^2e_1\bigr]
-s^2e_1^2e_2^2e_1-(rs)^3e_1e_2^2e_1^2\\
&=0.
\endsplit
$$
Thus, we get the identity (6)
$$\Cal T_1(e_2)\Cal T_1'(e_2)=r^2\Cal T_1'(e_2)\Cal
T_1(e_2).$$

To check the identity (7), we have
$$\split
e_1\Cal
T_2(e_1)^2&=e_1^2e_2e_1e_2-r^2e_1^2e_2^2e_1-r^2e_1e_2e_1^2e_2+r^4e_1e_2e_1e_2e_1,\\
-s(r{+}s)\Cal T_2(e_1)e_1\Cal T_2(e_1)&=
-s(r{+}s)e_1e_2e_1^2e_2+r^2s(r{+}s)e_1e_2e_1e_2e_1\\
&\quad+
r^2s(r{+}s)e_2e_1^3e_2-r^4s(r{+}s)e_2e_1^2e_2e_1,\\
rs^3\Cal
T_2(e_1)^2e_1&=rs^3e_1e_2e_1e_2e_1-(rs)^3e_1e_2^2e_1^2-(rs)^3e_2e_1^2e_2e_1+r^5s^3e_2e_1e_2e_1^2.
\endsplit
$$
So we obtain
$$\split
\text{LHS of (7)}&= e_1^2e_2e_1e_2-r^2e_1(e_1e_2^2)e_1-\Delta
(e_1e_2e_1^2)e_2+r(r{+}s)(r^2{+}s^2)e_1e_2e_1e_2e_1\\
&\quad+r^2s(r{+}s)e_2e_1^3e_2-(rs)^3e_1e_2^2e_1^2+r^5s^3e_2e_1e_2e_1^2-r^3s\Delta
e_2(e_1^2e_2e_1)\\
&=
e_1^2e_2e_1e_2-r^2\bigl[(r^2{+}s^2)e_1e_2e_1e_2e_1-(rs)^2e_1e_2^2e_1^2\bigr]\\
&\quad- \bigl[(rs)^2e_2e_1^3e_2+(rs)^{-1}\Delta
e_1^2e_2e_1e_2-(rs)^{-1}e_1^3e_2^2\bigr]\\
&\quad+r(r{+}s)(r^2{+}s^2)e_1e_2e_1e_2e_1+r^2s(r{+}s)e_2e_1^3e_2-(rs)^3e_1e_2^2e_1^2\\
&\quad+r^5s^3e_2e_1e_2e_1^2-r^3s\bigl[(rs)\Delta
e_2e_1e_2e_1^2+e_2e_1^3e_2-(rs)^3e_2^2e_1^3\bigr]\\
&=-(rs)^{-1}(r^2{+}s^2)e_1^2(e_2e_1e_2)+(rs)(r^2{+}s^2)e_1(e_2e_1e_2)e_1+r^3s^2(r{-}s)e_1e_2^2e_1^2\\
&\quad+(rs)^{-1}e_1^3e_2^2-r^4s^2(r^2{+}s^2)(e_2e_1e_2)e_1^2+r^6s^4e_2^2e_1^3\\
&=-(rs)^{-1}\bigl[(rs)^2e_1^2e_2^2e_1+e_1^3e_2^2\bigr]
+(rs)\bigl[(rs)^2e_1e_2^2e_1^2+e_1^2e_2^2e_1\bigr]\\
&\quad+r^3s^2(r{-}s)e_1e_2^2e_1^2+(rs)^{-1}e_1^3e_2^2
-r^4s^2\bigl[(rs)^2e_2^2e_1^3+e_1e_2^2e_1^2\bigr]+r^6s^4e_2^2e_1^3\\
&=0.
\endsplit
$$

This complete the proof.\hfill\qed
\enddemo

Now let us consider the rank $3$ cases.

\proclaim{Lemma 3.7} \ For the rank $3$ cases of types $A_3$,
$B_3$ and $C_3$, the Lusztig's symmetries $\Cal T_i$ preserve the
defining relations $(X2)$ or $(X3)$ of $(U_{r,s}(\frak g),
\langle\,,\rangle)$ into its associated object $(U_{s^{-1},
r^{-1}}(\frak g), \langle\,|\,\rangle)$, for $X=A,\, B,\,C$ if and
only if $rs=1$.
\endproclaim
\demo{Proof} \ When $rs=1$,  we set $r=q$, $s=q^{-1}$.
It is then obvious that $U_{q,q^{-1}}(\frak g)$ has Lusztig's
symmetries for any type and any rank of $\frak g$. Indeed, the quantum
group $U_q(\frak g)$ of Drinfel'd-Jimbo type is a quotient of
$U_{q,q^{-1}}(\frak g)$ by the ideal $(\om_i^{-1}-\om_i',
\;i=1,\cdots,n)$ (see \cite{15}). Consequently, each $\Cal T_i$ automatically preserves the
defining relations $(X2)$ and $(X3)$.

Conversely, assume that each $\Cal T_i$ preserves relations $(X2)$
or $(X3)$. In type $A_3$ or $C_3$ case, we consider the situation
where $\Cal T_3$ preserves $\om_1e_2\om_1^{-1}=s\,e_2$. This leads
to $\Cal T_3(\om_1)\Cal T_3(e_2)\Cal
T_3(\om_1^{-1})=\om_1(e_2e_3-re_3e_2)\om_1^{-1}=s\Cal
T_3(e_2)=\Cal T_3(\om_1e_2\om_1^{-1})=\Cal T_3(s\,e_2)=r^{-1}\Cal
T_3(e_2)$ and this implies $rs=1$. Similarly, in $B_3$ case, we
can get the same condition provided that we observe that $\Cal
T_1$ preserves the relation $\om_3e_2\om_3^{-1}=r^{-1}e_2$.
\hfill\qed
\enddemo

\medskip

\demo{Proof of Theorem 3.1} \ When $\text{rank}\,(\frak g)=2$,
Lemmas 3.3, 3.4 \& 3.5 indicate that for arbitrary parameters $r$,
$s$ with $r^2\ne s^2$, $U_{r,s}(\frak {sl}_3)$, $U_{r,s}(\frak
{sp}_4)$ and $U_{r,s}(\frak {so}_5)$ possess the Lusztig
symmetries into their respective associated quantum groups. When
$\text{rank}\,(\frak g)>2$, this $\frak g$ contains one of rank
$3$ Lie subalgebras $\frak{sl}_4$, $\frak{so}_7$ and $\frak
{sp}_6$. Lemma 3.7 gives the required assertion. \hfill\qed
\enddemo

\bigskip\bigskip
\heading{$\bold 4$. \ Appendix: Some
Calculations}\endheading
\bigskip

\noindent {\bf 4.1 \ Calculations in the proof of Proposition 2.3}
\ The relevant terms of $\Delta(X)$ for $X=f_1^3f_2$
 in (6C) are as
follows:
$$\split
&f_1{\om_1'}^2\om_2'\ot f_1\om_1'\om_2'\ot f_1\om_2'\ot
f_2+f_1{\om_1'}^2\om_2'\ot\om_1'f_1\om_2'\ot f_1\om_2'\ot f_2
\\
&\qquad+\om_1'f_1\om_1'\om_2'\ot f_1\om_1'\om_2'\ot f_1\om_2'\ot
f_2 +\om_1'f_1\om_1'\om_2'\ot\om_1'f_1\om_2'\ot f_1\om_2'\ot
f_2\\
&\qquad\qquad+{\om_1'}^2f_1\om_2'\ot f_1\om_1'\om_2'\ot
f_1\om_2'\ot f_2
+{\om_1'}^2f_1\om_2'\ot\om_1'f_1\om_2'\ot f_1\om_2'\ot f_2\\
& +f_1{\om_1'}^2\om_2'\ot f_1\om_1'\om_2'\ot\om_1'f_2\ot f_1
+f_1{\om_1'}^2\om_2'\ot\om_1'f_1\om_2'\ot\om_1'f_2\ot f_1
\\
&\qquad+\om_1'f_1\om_1'\om_2'\ot f_1\om_1'\om_2'\ot\om_1'f_2\ot
f_1 +\om_1'f_1\om_1'\om_2'\ot\om_1'f_1\om_2'\ot\om_1'f_2\ot
f_1\\
&\qquad\qquad+{\om_1'}^2f_1\om_2'\ot
f_1\om_1'\om_2'\ot\om_1'f_2\ot f_1
+{\om_1'}^2f_1\om_2'\ot\om_1'f_1\om_2'\ot\om_1'f_2\ot f_1\\
& +f_1{\om_1'}^2\om_2'\ot{\om_1'}^2f_2\ot f_1\om_1'\ot f_1
+f_1{\om_1'}^2\om_2'\ot{\om_1'}^2f_2\ot\om_1'f_1\ot f_1\\
&\qquad+\om_1'f_1\om_1'\om_2'\ot {\om_1'}^2f_2\ot f_1\om_1'\ot f_1
+ \om_1'f_1\om_1'\om_2'\ot{\om_1'}^2f_2\ot \om_1'f_1\ot f_1\\
&\qquad\qquad+{\om_1'}^2f_1\om_2'\ot{\om_1'}^2f_2\ot f_1\om_1'\ot
f_1+{\om_1'}^2f_1\om_2'\ot{\om_1'}^2f_2\ot\om_1'f_1\ot
f_1\\
& +{\om_1'}^3f_2\ot f_1{\om_1'}^2\ot f_1\om_1'\ot f_1
+{\om_1'}^3f_2\ot f_1{\om_1'}^2\ot\om_1'f_1\ot f_1\\
&\qquad+{\om_1'}^3f_2\ot {\om_1'}f_1\om_1'\ot f_1\om_1'\ot f_1
+{\om_1'}^3f_2\ot \om_1'f_1\om_1'\ot\om_1'f_1\ot f_1\\
&\qquad\qquad+{\om_1'}^3f_2\ot{\om_1'}^2f_1\ot f_1\om_1'\ot
f_1+{\om_1'}^3f_2\ot{\om_1'}^2f_1\ot\om_1'f_1\ot f_1.
\endsplit\tag{4.1}
$$

The relevant terms of $\Delta(X)$ for $X=f_2f_1^3$ in (6C) are as
follows:
$$
\split &\om_1'\om_2'f_1\om_1'\ot\om_1'\om_2'f_1\ot\om_2'f_1\ot
f_2+{\om_1'}^2\om_2'f_1\ot\om_1'\om_2'f_1\ot\om_2'f_1\ot f_2\\
&\qquad+\om_1'\om_2'f_1\om_1'\ot\om_2'f_1\om_1'\ot\om_2'f_1\ot f_2
+{\om_1'}^2\om_2'f_1\ot\om_2'f_1\om_1'\ot\om_2'f_1\ot f_2\\
&\qquad\qquad+\om_2'f_1{\om_1'}^2\ot\om_2'f_1\om_1'\ot\om_2'f_1\ot
f_2+\om_2'f_1{\om_1'}^2\ot\om_1'\om_2'f_1\ot\om_2'f_1\ot f_2\\
&+\om_1'\om_2'f_1\om_1'\ot\om_1'\om_2'f_1\ot f_2\om_1'\ot f_1
+{\om_1'}^2\om_2'f_1\ot\om_1'\om_2'f_1\ot f_2\om_1'\ot f_1\\
&\qquad+\om_1'\om_2'f_1\om_1'\ot\om_2'f_1\om_1'\ot f_2\om_1'\ot
f_1+{\om_1'}^2\om_2'f_1\ot\om_2'f_1\om_1'\ot f_2\om_1'\ot f_1\\
&\qquad\qquad+\om_2'f_1{\om_1'}^2\ot\om_2'f_1\om_1'\ot
f_2\om_1'\ot f_1+\om_2'f_1{\om_1'}^2\ot\om_1'\om_2'f_1\ot
f_2\om_1'\ot f_1\\
&+\om_2'f_1{\om_1'}^2\ot f_2{\om_1'}^2\ot f_1\om_1'\ot f_1
+\om_2'f_1{\om_1'}^2\ot f_2{\om_1'}^2\ot \om_1'f_1\ot f_1\\
&\qquad+\om_1'\om_2'f_1\om_1'\ot f_2{\om_1'}^2\ot f_1\om_1'\ot f_1
+{\om_1'}^2\om_2'f_1\ot f_2{\om_1'}^2\ot f_1\om_1'\ot f_1\\
&\qquad\qquad+\om_1'\om_2'f_1\om_1'\ot f_2{\om_1'}^2\ot
\om_1'f_1\ot f_1+{\om_1'}^2\om_2'f_1\ot f_2{\om_1'}^2\ot
\om_1'f_1\ot f_1\\
&+f_2{\om_1'}^3\ot f_1{\om_1'}^2\ot f_1\om_1'\ot f_1 +
f_2{\om_1'}^3\ot f_1{\om_1'}^2\ot \om_1'f_1\ot f_1\\
&\qquad+f_2{\om_1'}^3\ot \om_1'f_1\om_1'\ot f_1\om_1'\ot f_1
+f_2{\om_1'}^3\ot {\om_1'}^2f_1\ot f_1\om_1'\ot f_1\\
&\qquad\qquad+f_2{\om_1'}^3\ot \om_1'f_1\om_1'\ot \om_1'f_1\ot f_1
+f_2{\om_1'}^3\ot {\om_1'}^2f_1\ot \om_1'f_1\ot f_1.
\endsplit\tag{4.2}
$$

The relevant terms of $\Delta(X)$ for $X=f_1^2f_2f_1$ in (6C) are
as follows:
$$\split
&f_1{\om_1'}^2\om_2'\ot f_1\om_1'\om_2'\ot\om_2'f_1\ot f_2 +
\om_1'f_1\om_1'\om_2'\ot f_1\om_1'\om_2'\ot \om_2'f_1\ot
f_2\\
&\qquad+f_1{\om_1'}^2\om_2'\ot\om_1'\om_2'f_1\ot f_1\om_2'\ot f_2+
\om_1'f_1\om_1'\om_2'\ot\om_1'\om_2'f_1\ot f_1\om_2'\ot f_2\\
&\qquad\qquad+{\om_1'}^2\om_2'f_1\ot f_1\om_1'\om_2'\ot
f_1\om_2'\ot f_2+ {\om_1'}^2\om_2'f_1\ot \om_1'f_1\om_2'\ot
f_1\om_2'\ot f_2 \\
&+f_1{\om_1'}^2\om_2'\ot f_1\om_1'\om_2'\ot f_2\om_1'\ot f_1+
\om_1'f_1\om_1'\om_2'\ot f_1\om_1'\om_2'\ot f_2\om_1'\ot f_1\\
&\qquad+f_1{\om_1'}^2\om_2'\ot\om_1'\om_2'f_1\ot\om_1'f_2\ot f_1+
\om_1'f_1\om_1'\om_2'\ot\om_1'\om_2'f_1\ot\om_1'f_2\ot f_1\\
&\qquad\qquad+{\om_1'}^2\om_2'f_1\ot
f_1\om_1'\om_2'\ot\om_1'f_2\ot f_1+ {\om_1'}^2\om_2'f_1\ot
\om_1'f_1\om_2'\ot\om_1'f_2\ot f_1\\
&+f_1{\om_1'}^2\om_2'\ot\om_1'f_2\om_1'\ot f_1\om_1'\ot
f_1+f_1{\om_1'}^2\om_2'\ot\om_1'f_2\om_1'\ot \om_1'f_1\ot f_1\\
&\qquad+\om_1'f_1\om_1'\om_2'\ot\om_1'f_2\om_1'\ot f_1\om_1'\ot
f_1
+\om_1'f_1\om_1'\om_2'\ot\om_1'f_2\om_1'\ot \om_1'f_1\ot f_1\\
&\qquad\qquad+{\om_1'}^2\om_2'f_1\ot{\om_1'}^2f_2\ot f_1\om_1'\ot
f_1+{\om_1'}^2\om_2'f_1\ot{\om_1'}^2f_2\ot \om_1'f_1\ot f_1 \\
&+{\om_1'}^2f_2\om_1'\ot f_1{\om_1'}^2\ot f_1\om_1'\ot f_1+
{\om_1'}^2f_2\om_1'\ot \om_1'f_1\om_1'\ot f_1\om_1'\ot f_1\\
&\qquad+{\om_1'}^2f_2\om_1'\ot f_1{\om_1'}^2\ot \om_1'f_1\ot f_1+
{\om_1'}^2f_2\om_1'\ot \om_1'f_1\om_1'\ot \om_1'f_1\ot f_1\\
&\qquad\qquad+{\om_1'}^2f_2\om_1'\ot {\om_1'}^2f_1\ot f_1\om_1'\ot
f_1+ {\om_1'}^2f_2\om_1'\ot {\om_1'}^2f_1\ot \om_1'f_1\ot f_1.
\endsplit\tag{4.3}
$$

The relevant terms of $\Delta(X)$ for $X=f_1f_2f_1^2$ in (6C) are
as follows:
$$\split
&f_1{\om_1'}^2\om_2'\ot\om_2'f_1\om_1'\ot\om_2'f_1\ot f_2 +
f_1{\om_1'}^2\om_2'\ot\om_1'\om_2'f_1\ot\om_2'f_1\ot f_2\\
&\qquad+\om_1'\om_2'f_1\om_1'\ot f_1\om_1'\om_2'\ot\om_2'f_1\ot
f_2+ {\om_1'}^2\om_2'f_1\ot f_1\om_1'\om_2'\ot\om_2'f_1\ot f_2\\
&\qquad\qquad+\om_1'\om_2'f_1\om_1'\ot \om_1'\om_2'f_1\ot
f_1\om_2'\ot f_2+ {\om_1'}^2\om_2'f_1\ot \om_1'\om_2'f_1\ot
f_1\om_2'\ot f_2\\
&+f_1{\om_1'}^2\om_2'\ot\om_2'f_1\om_1'\ot f_2\om_1'\ot f_1+
f_1{\om_1'}^2\om_2'\ot\om_1'\om_2'f_1\ot f_2\om_1'\ot f_1\\
&\qquad+\om_1'\om_2'f_1\om_1'\ot f_1\om_1'\om_2'\ot f_2\om_1'\ot
f_1+ {\om_1'}^2\om_2'f_1\ot f_1\om_1'\om_2'\ot f_2\om_1'\ot f_1\\
&\qquad\qquad+\om_1'\om_2'f_1\om_1'\ot \om_1'\om_2'f_1\ot
\om_1'f_2\ot f_1+{\om_1'}^2\om_2'f_1\ot \om_1'\om_2'f_1\ot
\om_1'f_2\ot f_1\\
&+f_1{\om_1'}^2\om_2'\ot f_2{\om_1'}^2\ot f_1\om_1'\ot f_1+
f_1{\om_1'}^2\om_2'\ot f_2{\om_1'}^2\ot \om_1'f_1\ot f_1\\
&\qquad+\om_1'\om_2'f_1\om_1'\ot\om_1'f_2\om_1'\ot f_1\om_1'\ot
f_1+{\om_1'}^2\om_2'f_1\ot\om_1'f_2\om_1'\ot f_1\om_1'\ot f_1\\
&\qquad\qquad+\om_1'\om_2'f_1\om_1'\ot\om_1'f_2\om_1'\ot
\om_1'f_1\ot f_1+{\om_1'}^2\om_2'f_1\ot\om_1'f_2\om_1'\ot
\om_1'f_1\ot f_1\\
&+\om_1'f_2{\om_1'}^2\ot f_1{\om_1'}^2\ot f_1\om_1'\ot f_1+
\om_1'f_2{\om_1'}^2\ot f_1{\om_1'}^2\ot \om_1'f_1\ot f_1\\
&\qquad+ \om_1'f_2{\om_1'}^2\ot \om_1'f_1\om_1'\ot f_1\om_1'\ot
f_1+\om_1'f_2{\om_1'}^2\ot {\om_1'}^2f_1\ot f_1\om_1'\ot f_1\\
&\qquad\qquad+\om_1'f_2{\om_1'}^2\ot \om_1'f_1\om_1'\ot
\om_1'f_1\ot f_1+\om_1'f_2{\om_1'}^2\ot {\om_1'}^2f_1\ot
\om_1'f_1\ot f_1.
\endsplit\tag{4.4}
$$

The relevant terms of $\Delta^{(2)}(f_2^2f_1)$ in (7C) are as
follows:
$$\split
&f_2\om_1'\om_2'\ot f_2\om_1'\ot f_1+\om_2'f_2\om_1'\ot f_2\om_1'\ot f_1\\
&\qquad+f_2\om_1'\om_2'\ot \om_2'f_1\ot f_2+\om_2'f_2\om_1'\ot
\om_2'f_1\ot f_2\\
&\qquad\qquad+{\om_2'}^2f_1\ot f_2\om_2'\ot f_2+{\om_2'}^2f_1\ot
\om_2'f_2\ot f_2.
\endsplit\tag{4.5}
$$

The relevant terms of $\Delta^{(2)}(f_2f_1f_2)$ in (7C) are as
follows:
$$\split
&f_2\om_1'\om_2'\ot \om_1' f_2\ot f_1+\om_1'\om_2'f_2\ot
f_2\om_1' \ot f_1\\
&\qquad+f_2\om_1'\om_2'\ot f_1\om_2'\ot f_2+\om_1'\om_2'f_2\ot
\om_2'f_1\ot f_2\\
&\qquad\qquad+\om_2'f_1\om_2'\ot f_2\om_2'\ot
f_2+\om_2'f_1\om_2'\ot \om_2'f_2\ot f_2.
\endsplit\tag{4.6}
$$

\bigskip

\bigskip

\bigskip

\heading{\bf Acknowledgments}\endheading

\medskip

The study of the two-parameter quantum groups of types $B,\,C,\,D$
was early initiated when NH visited l'Institut de Recherche
Math\'ematique Avanc\'ee (IRMA -- C.N.R.S), Strasbourg in the
academic year 2000-2001, he would like to thank Prof. C. Kassel
for his hospitality. YG would like to thank Professors M.
Kashiwara and K. Saito for their extreme hospitality during his
visit to RIMS of Kyoto University when part of this work was done.
A special thanks go to Honglian Zhang for verifying Lemma 3.6.

\end